\documentclass[10pt]{elsarticle}
\usepackage{fancyhdr}
\usepackage{marginnote}

\usepackage{amsmath}
\usepackage{tikz}
\usepackage{frcursive}
\usepackage{mathtools}
\usepackage{amsfonts}
\usepackage{mathrsfs}
\usepackage{amsthm}
\usepackage{amssymb}
\usepackage{color}
\usepackage{dsfont}
\usepackage{geometry}
\newgeometry{tmargin=2.8cm, bmargin=3.5cm, lmargin=1.7cm, rmargin=1.7cm}

\usepackage{color}

\DeclareMathAlphabet{\mathpzc}{OT1}{pzc}{m}{it}

\newtheorem{theo}{\bf Theorem} 
 
\newtheorem{coro}{\bf Corollary}[section]
\newtheorem{lem}[coro]{\bf Lemma} 
\newtheorem{rem}[coro]{\bf Remark}

\newtheorem{prop}[coro]{\bf Proposition}

\newcommand{\blue}[1]{\textcolor{black}{#1}}
\newcommand{\opA}{{\mathcal{ A}}}

\def\avenorm#1{\mathchoice%
          {\mathop{\kern 0.2em\vrule width 0.6em height 0.69678ex depth -0.58065ex
                  \kern -0.545em \|{#1}\|}}%
          {\mathop{\kern 0.1em\vrule width 0.5em height 0.69678ex depth -0.60387ex
                  \kern -0.495em \|{#1}\|}}%
          {\mathop{\kern 0.1em\vrule width 0.5em height 0.69678ex depth -0.60387ex
                  \kern -0.495em \|{#1}\|}}%
          {\mathop{\kern 0.1em\vrule width 0.5em height 0.69678ex depth -0.60387ex
                  \kern -0.495em \|{#1}\|}}}

\newcommand{\barint}{
         \rule[.036in]{.12in}{.009in}\kern-.16in
          \displaystyle\int  } 
          
\newcommand{\data}{\textit{\texttt{data}}}
\def\R{{\mathbb{R}}}

\def\r{{\mathbb{R}}}
\def\N{{\mathbb{N}}}
\def\rn{{\mathbb{R}^{n}}}

\newcommand{\oT}{{\overline{\mathcal{T}}}}
\newcommand{\supp}{{\mathrm{supp}}}

\def\rp{{[0,\infty )}}

\def\ve{{\varepsilon}}
\def\vr{{\varrho}}

\def\cW{{\mathpzc{W}}}
\def\bW{{\mathbf{W}}}
\def\cI{{\mathcal{I}}}

\def\cV{{\mathpzc{V}}}
\def\bV{{\mathbf{V}}}
\def\cIa{{\cI_\alpha}}

\def\Waps{{\bW_{\alpha,\psi}}}
\def\Wap{{\cW_{\alpha,p}}}
\def\WaG{{\cW_{\alpha,G}}}
\def\WapsR{{\bW_{\alpha,\psi}^R}}

\def\WaGR{{\cW_{\alpha,G}^R}}

\def\xiL1{{\xi_{L^1}}}

\newcommand{\wt}{\widetilde}

\newcommand{\vp}{\varphi}

\newcommand{\dv}{\mathrm{div}}

\begin{document}

\begin{frontmatter}

\title{Boundedness of Wolff-type potentials and applications to PDEs 
}

\author[1]{Michał Borowski}\ead{m.borowski@mimuw.edu.pl}
\author[1]{Iwona Chlebicka\corref{mycorrespondingauthor}}
\cortext[mycorrespondingauthor]{Corresponding author}
\ead{i.chlebicka@mimuw.edu.pl}
\author[1]{Błażej Miasojedow}\ead{b.miasojedow@mimuw.edu.pl}

\fntext[myfootnote]{M.B. is supported by the Ministry of Science and Higher Education project Szko\l{}a Or\l{}\'o{}w, project number 500-D110-06-0465160. 
}

\address[1]{Institute of Applied Mathematics and Mechanics,
University of Warsaw, ul. Banacha 2, 02-097 Warsaw, Poland
}

\begin{abstract} 
We provide a short proof of a sharp rearrangement estimate for a generalized version of a potential of Wolff--Havin--Maz'ya type. As a consequence, we prove a reduction principle for that integral operators, that is, a characterization of
those rearrangement invariant spaces between which the potentials are bounded via a one-dimensional inequality of Hardy-type. 

Since the special case of the mentioned potential is known to control precisely very weak solutions to a broad class of quasilinear elliptic PDEs of non-standard growth, we infer the local regularity properties of the solutions in rearrangement invariant
spaces for prescribed classes of data.
\end{abstract}
\begin{keyword}  potential estimates\sep quasilinear elliptic PDEs\sep rearrangement\sep regularity\sep Wolff potential

\fntext[myfootnote]{The authors claim no conflict of interest.}
\fntext[myfootnote]{Data sharing not applicable to this article as no datasets were generated or analysed during
the current study.}
 
\MSC[2020] 31C15 (35A25, 35J62, 46E30)
\end{keyword}

\end{frontmatter}

\setcounter{tocdepth}{1}

\section{Introduction}  The classical Riesz potential and its nonlinear counterparts of Havin--Maz'ya and Wolff-type play a fundamental role in the theory of regularity of PDEs via the methods of harmonic analysis, see \cite{HeKiMa,KuMiguide} and references therein. We study properties of a broader class of operators than the mentioned above to give a tool for analysis of solutions to strongly nonlinear PDE problems and related minimizers to the calculus of variations. The corresponding methods of~potential analysis have recently found significant attention~\cite{KuMi2018jems,KuMiguide,CGZG-Wolff,CYZG-Wolff,CiSch,DM-Lip2021,Trudinger-Wang,Ed,DuMiJFA2010
,BeckMingione,Baroni-Riesz}. We prove a~general version of the classical reduction principle yielding the equivalence of the boundedness of nonlinear potentials to a~one-dimensional Hardy-type inequality. Main focus is on the class of potentials under minimal structural assumptions. 

We obtain the result resembling the bounds for the Riesz potential of \cite{ONeil} and for Havin--Maz'ya potential with a power $1<p<\infty$ of~\cite{Ci-pot}. Nonetheless, the potentials we consider do not share a structure involving Riesz potential as in~\cite{Ci-pot} and no growth condition of $\Delta_2$-type is imposed on the function used in the definition of the potential. Thus, it is impossible to employ the arguments of~\cite{Ci-pot} basing on the representation formulas for $K$-functionals from~\cite{holmstedt} (or~\cite{BeSa88}). To be more precise, we provide a short novel proof of a sharp rearrangement estimate for a general version of a potential of Wolff--Havin--Maz'ya-type generated by any non-decreasing and left-continuous function $\psi : \rp \to \rp$, see Theorem~\ref{theo:boundedness-psi}. As a consequence, we infer characterizations of those rearrangement invariant spaces between which the generalized potentials are bounded, see Theorems~\ref{theo:equivalence} and~\ref{theo:equivalence-operators}. Let us stress that to our best knowledge there are no available results on boundedness of nonlinear potentials of non-power growth even between the classical spaces. 

To stress applicability of the main accomplishments, in Section~\ref{sec:appl} we restrict attention to potentials with doubling growth and show the precise transfer of regularity from data to solutions to the Orlicz growth measure data problems including the scalar isotropic doubling case and vectorial doubling case.  Gradient regularity is also commented~on. Regularity to the problems with non-doubling growth is exposed in Section~\ref{sec:non-doub}. \newline

A well-established role in the analysis of the $p$-Laplace scalar and vectorial problems \cite{DuMiJFA2010,KiMa92,KuMi2018jems,KuMiguide} of the form \begin{equation}\label{plap}
    -\dv \big(|Du|^{p-2}Du\big)=f\,,\qquad 1<p<\infty\,,
\end{equation} is played by the Wolff potential defined as
\begin{equation}
    \label{wolff-potential-p}
 \Wap f(x):=\int_0^\infty \left(r^{\alpha p-n}{\int_{B(x,r)}|f(y)|\,dy}\right)^{\frac{1}{p-1}}\,\frac{dr}{r}= \int_0^\infty r^{\alpha-1}\left(r^{\alpha-n}\int_{B(x,r)}|f(y)|\,dy\right)^{\frac{1}{p-1}}\,dr\,.\end{equation}
 Note that for $p=2$ potential $\cW_{\alpha/2,p}$ is nothing but the classical Riesz potential
\[\cIa f(x):=\int_0^\infty r^{\alpha-n-1}{\int_{B(x,r)}|f{(y)}|\,dy}\,dr= c(\alpha,n) \int_\rn \frac{|f(x)|}{|x-y|^{n-\alpha}}\,dy\,,\] which is a standard tool in the analysis of the Poisson equation.  Potentials having similar properties to $\Wap$ and given by\begin{equation}
    \label{havin-mazya-potential}\cV_{\alpha,p}f(x):= \cIa\big(\cIa f \big)^{\frac{1}{p-1}} (x)\,,
\end{equation} were studied earlier by Havin and Maz'ya, see \cite{HaMa72}. Boundedness of $\cIa$ and $\cV_{\alpha,p}$ between rearrangement invariant spaces is studied e.g., in~\cite{Ci-pot}, later on the boundedness of $\cIa$ was revisited in~\cite{Ed} and of $\cV_{\alpha,p}$ in~\cite{CiSch}.

The available pointwise estimates for solutions to problems with Orlicz growth~\cite{CGZG-Wolff} imply not only that $\WaG$ given by
\begin{equation}\label{wolff-potential}
\WaG f(x):=\int_0^\infty r^{\alpha - 1}g^{-1}\left(r^{\alpha-n}{\int_{B(x,r)}|f(y)|\,dy}\right)\,dr\,
\end{equation}
generalizes the classical operators. It provides a precise relation between data and solution to a quasilinear elliptic problem \begin{equation}\label{eq-Or}
    -\dv \opA(x,Du)=f\,,
\end{equation}
 where the growth of $\opA$ with respect to the second variable is described by the means of a given $N$-function $G\in C^1 {((0,\infty))} \cap C([0, \infty))$ for which $g:=G'$ and $G\in\Delta_2\cap\nabla_2$. It is known since~\cite{Ta0,Ta1} that the theory of rearrangements is a very effective way in study of PDEs of Orlicz growth. We use an approach initiated in~\cite{Ci-pot} and study regularity of solutions via rearrangement estimate to a potential of the Wolff-type. There is a~pointwise estimate from below and from above of solution to~\eqref{eq-Or} involving the following version of the truncated potential $\WaGR f(x) := \int_{0}^{R}r^{\alpha - 1} g^{-1}\left(r^{\alpha-n}{\int_{B(x,r)}|f(y)|\,dy}\right)\,dr $. In fact, by \cite[Theorem~2]{CGZG-Wolff}  it is known that an $\opA$-superharmonic function $u$ generated by a nonnegative $f\in L^1_{loc}$ satisfies  
\begin{equation}
\label{est-wolff-G}
C_L\left(\cW_{1,G}^Rf (x)-R\right)\leq u(x )\leq C_U\left( \inf_{B(x,R)} u(x)+\cW_{1,G}^Rf (x)+R\right)
\end{equation}
with some $C_L,C_U>0$ depending only on the parameters of the problem and for all sufficiently small $R$. This result extends to the Orlicz setting the pre-eminent result for the $p$-Laplacian by Kilpel\"ainen and Mal\'y~\cite{KiMa92} that was studied later, e.g., in~\cite{Trudinger-Wang}. Note that this result is sharp in the sense that there does not exist any essentially better potential to substitute in place of $\cW_{1,G}^R$ as the same potential controls the solution from below and from above. Later on, these kinds of results were used to infer various local properties of solutions, see~\cite{KuMiguide} and~\cite{Ci-pot,CiSch,CGZG-Wolff}.
Having~\eqref{est-wolff-G}, it is visible that precise results on the boundedness of the potential directly implies information on the local regularity of solutions in Section~\ref{sec:regularity}.\newline

Let us summarize results formulated in Section~\ref{sec:theorems-psi} and proven in Section~\ref{sec:main-proofs}. With this aim, we recall the definitions of decreasing rearrangement and maximal rearrangement of a function due to~\cite{EKRea,BeSa88}. We define the decreasing rearrangement $f^\ast : [0, \infty) \to [0, \infty]$ of a measurable function $f:\rn\to\R$ by
\begin{equation}
    \label{def:rear}
f^\ast(t) := \sup\Big \{ s \geq 0 \colon \big|\{x\in \R^n:|f(x)|>s\}\big| > t  \Big\}\,,
\end{equation}
where we assume that $\sup\emptyset=0$. Equivalently it can also be defined as the function, which is right-continuous, non-increasing, and equimeasurable with $f$, i.e., $|\{x : |f(x)| > t\}| = |\{x: f^{*}(x) > t\}|$ for all $t > 0$. 
 The maximal rearrangement is defined by
\begin{equation}
    \label{def:max-rear}
f^{\ast \ast}(t): =\frac 1t \int_0^t f^\ast(s) \, ds\quad\text{and}\quad 
f^{\ast \ast}(0)=
f^{\ast}(0)\,.
\end{equation} We shall introduce a potential subordinated to a given non-decreasing and left-continuous function $\psi:\rp\to\rp$ which is not assumed to fall into the realm of the doubling regime. Let us define
\begin{equation}\label{wolff-potential-psi}
\Waps f(x):=\int_0^\infty r^{\alpha - 1}\psi\left(r^{\alpha-n}{\int_{B(x,r)}|f(y)|\,dy}\right)dr\,.
\end{equation}
Therefore for $\psi(s)=s^\frac{1}{p-1}$ we have $\Waps=\Wap$ defined in~\eqref{wolff-potential-p}, as well as for $G\in\Delta_2\cap\nabla_2$, $g=G'$, and $\psi(s)= g^{-1}(s)$, we have $\Waps=\WaG$ defined in~\eqref{wolff-potential}. 
We shall start with a rearrangement inequality satisfied by $\Waps$. It is known since~\cite{ONeil} that in the special case potential $\bW_{\alpha/2,{\rm id}}=\cW_{\alpha/2,2}=\cIa$ satisfies
\begin{equation}\label{eq:ONeil}(\cI_{\alpha}f)^*(t)\leq C\int_t^\infty s^{\frac{\alpha}{n}-1}f^{**}(s)\,ds\,,\quad\text{for some } \ C=C(\alpha,n) > 0\,.\end{equation}
This result was extended in~\cite[Theorem~2.1]{Ci-pot} to an estimate on nonlinear potential $\cV_{\alpha,p}$ for $1<p<\infty$ reading
\begin{equation}\label{eq:Cianchi} (\cV_{\alpha,p})^*(t)\leq C\int_{kt}^\infty s^{\frac{\alpha}{n}-1}\left(s^\frac{\alpha}{n} f^{**}(s)\right)^\frac{1}{p-1}\,ds\qquad\text{for }\ 0<\alpha p<n\end{equation}
and some $C,k>0$. We provide a similar result in Theorem~\ref{theo:boundedness-psi} for $\Waps$ under no conditions of doubling type imposed on the non-decreasing function $\psi$. Namely, we prove that there exist $C_{\bW,1},C_{\bW,2}>0$ such that
\[(\Waps f)^{*}(t) \leq C_{\bW,1} \int_{{t}}^{\infty} s^{\frac{\alpha}{n}-1}\psi\left(C_{\bW,2}s^{\frac{\alpha}{n}}f^{**}(s)\right)ds\,.\] 
The method of the proof of Theorem~\ref{theo:boundedness-psi} is shorter and more elementary than the proof of~\cite[Theorem~2.1]{Ci-pot}. In turn, it provides a way to simplify the reasoning in the classical power case. Let us note that generalised potentials $\bV_{\alpha,\psi}$ and $\Waps$ (as well as $\cV_{\alpha,p}$ and $\cW_{\alpha,p}$) are closely related, see Section~\ref{ssec:rem:V=W} for more comments on this issue. 

In the power growth case, the boundedness of potential $\cV_{\alpha,p}$ (given by~\eqref{havin-mazya-potential}) between rearrangement invariant spaces is known to be equivalent to a one-dimensional inequality for nonlinear Hardy-type operator, see~\cite{Ci-pot}. We present two results in this spirit. Theorem~\ref{theo:equivalence} is more natural in formulation, as it deals with norm boundedness  of nonlinear operators $\Waps:X\to Y$. On the other hand, Theorem~\ref{theo:equivalence-operators} provides a more general characterisation of boundedness of more natural application. See Remark~\ref{rem:exTx} for an illustration. We present a few particular consequences on boundedness of $\Waps$ in the end of Section~\ref{sec:main-proofs}. The ideas of the proofs are summarized later. 

Section~\ref{sec:appl} is devoted to illustration of the main results. Namely, we restrict attention to the boundedness of potentials given by doubling functions and investigate their applications. Study on the boundedness of $\WaG$ is presented in Section~\ref{sec:WaG}, whereas its consequences for regularity for solutions to PDE problems of a type~\eqref{eq-Or} are presented in Section~\ref{sec:regularity}. Estimates for a solution and for its gradient follow from a combination of the bounds for $\WaG$ with known potential estimates from~\cite{CGZG-Wolff} and~\cite{Baroni-Riesz}, and embracing the notable results of~\cite{KiMa92,KuMiARMA2013}. As a consequence, the local integrability theory for solutions and their
gradients, is reduced to one-dimensional Hardy-type inequalities. 

\section{Boundedness of potential $\Waps$ } \label{sec:theorems-psi}

We start with presenting the following sharp result on boundedness of $\Waps$. This result relates to~\cite[Theorem~2.1]{Ci-pot}. Let us stress that our attention is restricted to operators defined for measurable functions $f:\rn\to\r$ such that \begin{equation}
    \label{skoncz-wart}
    |\{x : | f (x)| > t\}| < \infty\qquad\text{for }\ t > 0\,.
\end{equation}
\begin{theo} \label{theo:boundedness-psi}
Suppose $\psi:\rp\to\rp$ is non-decreasing and left-continuous. Let $n\geq 1$ and $\alpha\in(0,n)$. Then there exist  $C_{\bW,1}=C_{\bW,1}(\alpha,n)>0$ and $C_{\bW,2}=C_{\bW,2}(\alpha,n)>0$ such that
\begin{equation}\label{eq:West-psi}(\Waps f)^{*}(t) \leq C_{\bW,1} \int_{{t}}^{\infty} s^{\frac{\alpha}{n}-1}\psi\left(C_{\bW,2}s^{\frac{\alpha}{n}}f^{**}(s)\right)ds\end{equation} whenever $f\in L^1_{loc}(\rn)$ is such that~\eqref{skoncz-wart} holds.\\
The result is sharp, in the sense that there exist a constants $C_{\bW, 3} = C_{\bW, 3}(\alpha, n) > 0$ and $C_{\bW, 4}(\alpha, n) > 0$ such that
\begin{equation}\label{eq:Wsharp-psi}
C_{\bW, 3}\int_{{t}}^{\infty} s^{\frac{\alpha}{n}-1}\psi\left(C_{\bW, 4} s^{\frac{\alpha}{n}}f^{**}(s)\right)ds \leq (\Waps f)^*(t)\,
\end{equation}
if  $f\in L^1_{loc}(\rn)$ is any nonnegative and radially decreasing function satisfying~\eqref{skoncz-wart}. 
\end{theo}
Previous proof of the above fact, provided for $\psi(t)=t^{\frac{1}{p-1}}$ in~\cite{Ci-pot}, was based on the interpolation theory and applied $K$-functionals, see~\cite{BeSa88,holmstedt}. Related facts about $K$-functionals are studied since~\cite{kree,peetre,holmstedt} and were applied to establish many results concerning boundedness of various operators, e.g., in~\cite{CiKePi,Kerman-Pick,Ci-pot}. Nonetheless, operators covered by the methods involving $K$-functionals need to have certain structure. Let us indicate the difference between the structure of $\cV_{\alpha,p}$ given by~\eqref{havin-mazya-potential} and $\cW_{\alpha,p}$ from \eqref{wolff-potential-p}. The lack of the structure is even more striking in the case of $\Waps$ involving a general non-decreasing and left-continuous function $\psi$. Thus, as much as the type of result is similar to those of~\cite{Ci-pot}, our proofs cannot be based on the arguments from the mentioned contributions. Our proof is novel and more straightforward even in the case of the classical Wolff potential given in~\eqref{wolff-potential-p}, cf. \cite{Ci-pot}. Main steps make use of the properties of the fractional maximal operator.
\begin{rem}\label{rem:Wpsifinite}In the power case, in \cite[Theorem 2.1]{Ci-pot}, there is an assumption $\alpha p<n$. In the case of \blue{$\alpha p\geq n$} inequality~\eqref{eq:Cianchi} takes form $\infty\leq \infty$.
Note that if there exist $x\in\rn$ and $r>0$ such that $\int_{B(x,r)}|f(y)|\,dy>0$ and 
\begin{equation*}
    \int^{\infty} \psi(t^{1 - \frac{n}{\alpha}})dt = \infty\,,
\end{equation*}
then for every $t\geq 0$ it holds $(\Waps f)^{*}(t)=\infty$ and the right-hand side of~\eqref{eq:West-psi} is also infinite. Thus, the result is true for $\alpha\in(0,n)$, but meaningful only provided the mentioned integral is convergent.
\end{rem}
We apply Theorem~\ref{theo:boundedness-psi} in the proof of our next accomplishment, that is, in establishing an equivalent condition to the boundedness of the operator $\Waps:X\to Y$ between arbitrary rearrangement invariant quasi-normed function spaces. The proof follows the ideas of reasoning of \cite[Theorem~2.3]{Ci-pot} adapted to the operators involving general function $\psi$.

We refer the reader to \cite{BeSa88} for a~comprehensive treatment of a framework of rearrangement invariant spaces. Let us briefly recall definitions needed in the upcoming theorem. 

A quasi-normed {Banach} function space $X (\Omega)$ on a measurable subset $\Omega$ of $\rn$ is a linear space of measurable functions on $\Omega$ equipped with a functional $\|\cdot\|_{X(\Omega)}$, a quasi-norm, enjoying the following properties:
\begin{itemize}
    \item[(i)]  $\|f\|_{ X (\Omega)} > 0$ if $f \neq 0$;\\
            $\|\lambda f\|_{ X (\Omega)} = |\lambda|\,\| f\|_{X (\Omega)}$ for every $\lambda\in\R$ and $f\in X (\Omega)$;\\
            $\|f + g\|_{X (\Omega)}\leq c(\| f\|_{ X (\Omega)}  +\| g\|_{ X (\Omega)}  )$ for some constant $c\geq 1$ and for every $f, g \in X (\Omega)$;
    \item[(ii)] $0 \leq |g|\leq | f |$ a.e. in $\Omega$ implies $\|g\|_{ X (\Omega)} \leq \|f\|_{ X (\Omega)} $;
    \item[(iii)] $0\leq f_k\nearrow f$ a.e. implies $\|f_k\|_{ X (\Omega)} \nearrow\|f\|_{ X (\Omega)} $ as $k\to\infty$;
    \item[(iv)] if $E$ is a measurable subset of $\Omega$ and $|E| < \infty$, then $\|\mathds{1}_E\|_{ X (\Omega)} < \infty$;
    \item[(v)] for every measurable $E\subset\Omega$ with $|E|<\infty$, there exists a constant $C$ such that \mbox{$\int_E | f |\,d x\leq C\| f\|_{ X (\Omega)}$} for every $f \in X (\Omega)$.
\end{itemize}
The space $X (\Omega)$ is called a~Banach function space if the triangle inequality holds true, i.e., in (i) we have $c = 1$. In this case, the functional $\|\cdot\|_{X(\Omega)}$ is actually a norm which makes $X (\Omega)$ a~Banach
space. 

We say that  $f, g$ are equimeasurable if $|\{x : |f(x)| > t\}| = |\{x : |g(x)| > t\}|$ holds for all $t > 0$. A~quasi-normed function space   $X (\Omega)$ is called rearrangement invariant if for every pair $f,g$ of functions being equimeasurable and finite a.e. it holds that $||f||_{X(\Omega)} = ||g||_{X(\Omega)}$. By Luxemburg Representation  Theorem~\cite[Chapter~2, Theorem~4.10]{BeSa88}, this condition implies the existence of rearrangement invariant space $\overline{X} (0, |\Omega|)$ having the property that
\begin{equation}
    \label{repsp} \|f\|_{X(\Omega)}=\|f^*\|_{\overline{X}(0,|\Omega|)}\qquad\text{for every } \ f\in X(\Omega)\,.
\end{equation}
Space $\overline{X}$ is called the representation space of $X$. Classical instances of rearrangement-invariant quasi-normed spaces are Lebesgue
spaces, Lorentz spaces, Lorentz--Zygmund spaces, Orlicz spaces.

Here, and in what follows, any function $f$ on $\Omega$ is understood to be continued by $0$ outside $\Omega$. Observe that if $X (\Omega)$ is a rearrangement invariant quasi-normed space, then
$\|f\|_{X (\Omega)}=\|g\|_{X (\Omega)}$ if $f^* = g^*$. \newline

We provide the following result on characterisation of rearrangement invariant spaces between which $\Waps$ is bounded. As a matter of fact, the boundedness of $\Waps:X\to Y$  is proven to be equivalent to a one-dimensional Hardy-type inequality in the representation spaces of  $X$ and $Y$.

\begin{theo}\label{theo:equivalence} Suppose $\psi:\rp\to\rp$ is a non-decreasing and left-continuous function, $n\geq 1$, $\alpha\in(0,n)$ and let $\Omega \subseteq \rn$ be a measurable set. Let $X(\Omega)$ and $Y(\Omega)$ be quasi-normed Banach rearrangement invariant spaces. Then from the fact that
\begin{enumerate}
   \item[(i)] there exists a constant $c>0$ such that for every nonnegative function $\phi \in \overline{X}(0,|\Omega|)$ it holds
    \begin{equation*} \left\| \int_{{t}}^{\infty} s^{\frac{\alpha}{n}-1}\psi\left(s^{\frac{\alpha}{n}-1}\int_{0}^{s}\phi(y)\,dy\right)ds \right\|_{\overline{Y}(0,|\Omega|)} \leq c ||\phi||_{\overline{X}(0,|\Omega|)}\,, \end{equation*}
\end{enumerate}
it follows that
\begin{enumerate}
    \item[(ii)] there exists a constant $c>0$ such that for every $f \in \blue{X(\Omega)}$ it holds that
    \[||\Waps f||_{Y(\Omega)} \leq c ||f||_{X(\Omega)}\,.\]
\end{enumerate} 
Moreover, for $\Omega = \rn$ also the reverse implication holds true.
\end{theo}
\blue{In Section~\ref{sec:WaG} we provide the counterpart of the equivalence  {\it (ii)$\iff $(i)}  when $|\Omega| < \infty$ and the function defining the potential satisfies doubling growth conditions, see Theorem~\ref{theo:equivalence-G}.}

Because of the structural similarity between Theorems~\ref{theo:boundedness-psi} and~\ref{theo:equivalence} and their counterparts in~\cite{Ci-pot}, in the proof of Theorem~\ref{theo:equivalence} it was possible to make use of the arguments of \cite[Theorem~2.3]{Ci-pot}. Let us point out that the third equivalent item of this theorem does not follow because of the different structural properties of the potential considered there and of $\Waps$ with $\psi$ of non-power growth.

Note however that, unlike Lebesgue and Lorentz spaces, where using a norm is typical, in the Orlicz and Lorentz--Orlicz-type settings, it is more natural to handle modulars rather than norms. Therefore we extend Theorem~\ref{theo:equivalence-operators} to the following form using similar yet  arguments.

\begin{theo}\label{theo:equivalence-operators} Suppose $\psi:\rp\to\rp$ is a non-decreasing and left-continuous function, $n\geq 1$, and $\alpha\in(0,n)$. \blue{Let $X(\rn)$ and $Y(\rn)$} be quasi-normed rearrangement invariant spaces. Suppose \blue{$\mathcal{T}_Y:Y(\rn)\to\rp$,  $\mathcal{T}_X:X(\rn)\to\rp$,  $\oT_{\overline{Y}}:\overline{Y}(0,\infty)\to\rp$, $\oT_{\overline{X}}:\overline{X}(0,\infty)\to\rp$} satisfy\begin{equation}\label{oT-def}
    \mathcal{T}_X(f_1)=\oT_{\overline{X}}(f_1^*)\ \text{ for }f_1\in X(\blue{\rn})  \qquad\text{and}\qquad  \mathcal{T}_Y(f_2)=\oT_{\overline{Y}}(f_2^*)\ \text{ for }f_2\in Y(\blue{\rn}) \,.
\end{equation}
Assume further that\begin{equation}
    \label{oT-rear} \text{if $\quad 0\leq \phi\in \overline{X}\blue{(0,\infty)},\quad$ then $\quad\oT_{\overline{X}}(\phi) = \oT_{\overline{X}}( \phi^*)$ \quad and \quad  $\quad\oT_{\overline{Y}}(\phi) = \oT_{\overline{Y}}( \phi^*)\,$.}
\end{equation} Additionally, suppose that
\begin{equation}
    \label{oTY-monot} \text{if $\quad \phi_1,\phi_2\in \overline{Y}\blue{(0,\infty)}\quad$ and $\quad 0\leq \phi_1\leq \phi_2,\quad$ then $\quad\oT_{\overline{Y}}(\phi_1)\leq  \oT_{\overline{Y}}(\phi_2)\,$.}
\end{equation}
\blue{Let $\Omega \subseteq \rn$ be a measurable set. } Then from  the fact that
\begin{enumerate}
    \item[(i)] for every $k > 0$ there exists a constant $c>0$ such that for every nonnegative function $\phi \in \overline{X}(0, |\Omega|)$ it holds
    \begin{equation*} \oT_{\overline{Y}}\left(k\int_{{t}}^{\infty} s^{\frac{\alpha}{n}-1}\psi\left(s^{\frac{\alpha}{n}-1}\int_{0}^{s}\phi(y)dy\right)ds\right) \leq \oT_{\overline{X}}(c\phi)\,, \end{equation*}   \end{enumerate}
it follows that
\begin{enumerate}
    \item[(ii)] for every $k>0$ there exists a constant $c > 0$ such that for every $f \in \blue{X(\Omega)}$ it holds that
    \[\mathcal{T}_Y(k\Waps f) \leq  \mathcal{T}_X(cf)\,.\]
   \end{enumerate}
Moreover, for $\Omega = \rn$ also the reverse implication holds true.
\end{theo}

\begin{rem}\label{rem:exTx}
Admissible choice of $\mathcal{T}_X,\mathcal{T}_Y$ in Theorem~\ref{theo:equivalence-operators} is  \[\mathcal{T}_Xf=\int_0^\infty s^\sigma H_X(s^\vr f^*(s))\,ds\quad\text{ and }\quad\mathcal{T}_Yf=\int_0^\infty s^\alpha H_Y(s^\beta f^*(s))\,ds\] for some non-decreasing functions $H_X,H_Y:\rp\to\rp$ and $\alpha,\beta, \sigma, \vr \in \R$. Indeed, let us fix $\phi\in\overline{X}$ and define  $f_\phi(x):=\phi(\omega_n|x|^n)$. Then $f^*_\phi(s)=\phi^*(s)$. We denote $\oT_{\overline{X}}\phi:=\mathcal{T}_Xf_\phi.$ Then
\[\oT_{\overline{X}}f^*=\mathcal{T}_Xf_{f^*}=\int_0^\infty s^\sigma H_X(s^\vr f^*(s))\,ds\,,\]
therefore~\eqref{oT-def} holds true. On the other hand, to justify~\eqref{oT-rear} we notice that
\[\oT_{\overline{X}}\phi=\mathcal{T}_Xf_\phi=\int_0^\infty s^\sigma H_X(s^\vr (f_\phi)^*(s))\,ds=\int_0^\infty s^\sigma H_X(s^\vr \phi^*(s))\,ds\,.\]
Moreover,~\eqref{oTY-monot} follows from the monotonicity of $H_Y$.
\end{rem}

\section{Preliminaries}

Following a usual custom, we denote by $c$ or $C$ a general positive constant. Different occurrences from line to line will be still denoted by $c$ or $C$. Special occurrences will be denoted by $c_1, C_2$ or similarly. In order to stress relevant dependencies on parameters, we will use parentheses, i.e., $c= c(p,\alpha)$ means that $c$ depends on $p,\alpha$. Moreover, we say that two functions $f_1, f_2 : [0, \infty) \to \R$ are comparable and denote it as $f_1(t) \approx f_2(t)$ if there exist constants $c_1, c_2 > 0$ such that $f_1(c_1t) \leq f_2(t) \leq f_1(c_2t)$ for all $t > 0$. 
 
Let us denote the characteristic function of the set $E$ by $\mathds{1}_E$. We denote by $B(x,r)$ a ball centred in $x$ and radius $r>0$. By $\omega_n$ we denote a volume of a unit ball, i.e., $\omega_n:=|B(0,1)|$. 
  
  Given a measurable subset $E$ of $\Omega$, we denote   $\|f\|_{X (E)} = \|f\mathds{1}_E\|_{X (\Omega)}$
for any measurable function $f$ on $\Omega$. Note that then $||f||_{X(E)} \neq ||f^*||_{\overline{X}(0, |E|)}$ but rather $||f||_{X(E)} = ||\left(f \mathds{1}_{E}\right)^*||_{\overline{X}(0, |E|)}$. 
Moreover, we denote by $X_{loc} (\Omega)$ the space of measurable functions $f$ in $\Omega$
such that $\|f\|_{X (E)}<\infty$ for every compact set $E\subset\Omega$.

We use the basic properties of the decreasing rearrangement and the maximal rearrangement, defined in~\eqref{def:rear} and~\eqref{def:max-rear}. While considering these objects, it is usually assumed that~\eqref{skoncz-wart}~holds. Observe that 
\begin{equation}\label{eq:monrea} 
    \text{if $f, g : \rn \to \R$ and $f(x) \leq g(x)$ for all $x \in \rn,\ \ \ $ then } f^*(t) \leq g^*(t) \text{ for all $t > 0$}\,.
\end{equation}
By \cite[Chapter 1, Proposition 1.7]{BeSa88}, for any $f_1,f_2$ and $t>0$ it holds that
\begin{equation}\label{rearar-sum}
(f_1 + f_2)^*(t) \leq f_1^*(\tfrac{t}{2}) + f_2^*(\tfrac{t}{2}) \,.
\end{equation} The function $f^{**}$ is nonnegative, non-increasing and the function $t \mapsto tf^{**}(t)$ is non-decreasing. In particular
\begin{equation}
    \label{max-rear-prop}
    \tfrac{1}{2}f^{**}(t) \leq f^{**}(2t) \leq f^{**}(t)\,.
\end{equation}
Moreover, $f^* \leq f^{**}$. For any measurable set $\Omega \subset \R^n$ the following inequality holds true
\begin{equation}\label{eq:decinq}
    \int_{\Omega} |f(y)|dy \leq \int_{0}^{|\Omega|} f^{*}(s)ds = |\Omega|f^{**}(|\Omega|)\,.
\end{equation}
If $f$ is nonnegative, radially decreasing function, then 
\begin{equation}\label{eq:radec}f(x) = f^*(\omega_n |x|^n)\,\end{equation}
and 
\begin{equation}\label{eq:deceq}
    \int_{B(0, r)} f(y)dy = \int_{0}^{\omega_nr^n} f^{*}(s)ds = \omega_n r^nf^{**}(\omega_n r^n)\qquad \text{for every $r > 0\,$.}
\end{equation}
Lemma~\ref{lem:psi-rearr} gives us that for any non-decreasing, left-continuous function $\psi : [0, \infty) \to [0 ,\infty)$ we have the following identity
\begin{equation}\label{eq:psi-rearr}
    \left( \psi(|f|)\right)^{*}(t) = \psi(f^{*}(t))\,.
\end{equation} 


We use the classical fractional maximal operator $M_\alpha$, defined for every $f \in L^1_{loc}(\mathbb{R}^n)$ as
\begin{equation}\label{eq:maxop}
M_{\alpha}f(x) := \sup_{B \ni x} |B|^{\frac{\alpha}{n} - 1}\int_{B} |f(y)|dy\,,
\end{equation}
where $\alpha \in [0, n)$ is a parameter and supremum is taken over all balls $B$ such that $x \in B$, see~\cite{Stein}.
The following result has been proved in~\cite[Theorem~1.1]{Cimaximal}. For $n \geq 1$ and $\alpha \in [0, n)$, there \blue{exists} a constant $C_M = C_M(\alpha, n)$ such that for every $f \in L^1_{loc}(\mathbb{R}^n)$ it holds that
\begin{align}
    \label{remax} 
(M_{\alpha}f)^{*}(t) \leq C_M\sup_{s > t}s^{\frac{\alpha}{n}}f^{**}(s).
\end{align}
The proof of this fact employs basic properties of maximal function and rearrangements.


We use the Hardy inequality~\cite[Lemma 3.9]{BeSa88} of the following form. If $q \geq 1$ and $p < q - 1$, then for every measurable, nonnegative function $\phi : [0, \infty) \to [0, \infty)$ it holds that
\begin{equation}\label{eq:Hardy1}
    \int_{0}^{\infty} t^p \left( \frac{1}{t}\int_{0}^{t} \phi(x)\,dx \right)^q\,dt \leq \left(\frac{q}{q-p-1}\right)^q\int_{0}^{\infty}t^p\phi(t)^q\,dt\,.
\end{equation}
If instead of $p < q - 1$ we have $p > q - 1$, then
\begin{equation}\label{eq:Hardy2}
    \int_{0}^{\infty} t^p \left( \frac{1}{t}\int_{t}^{\infty} \phi(x)\,dx \right)^q\,dt \leq \left(\frac{q}{p + 1 - q}\right)^q\int_{0}^{\infty}t^p\phi(t)^q\,dt\,.
\end{equation}

Let $0 < p < \infty, 0 < q \leq \infty$. The Lorentz space $\Lambda^{p, q}(\Omega)$ is a space of measurable functions $f : \Omega \to \R$ such that
\begin{equation}\label{lorentznorm}
\|f\|_{{\Lambda^{p,q}}(\Omega)}:= \big\|s^{\frac 1p- \frac 1q}f^*(s)\big\|_{L^q(0, |\Omega|)} < \infty\,.
\end{equation}
 For $q=\infty$ this space is usually called the Marcinkiewicz space. Furthermore,  $\Lambda^{[p,q]}(\Omega)$ for $0 < p < \infty$, $0 < q \leq \infty$, denotes the Lorentz space equipped with the quasi-norm given by
\begin{equation}\label{lorentznorm**}
\|f\|_{\Lambda^{[p,q]}(\Omega)}:=\big\|s^{\frac 1p- \frac 1q}f^{**}(s)\big\|_{L^q (0, |\Omega|)}
\end{equation}
 for measurable functions $f : \Omega \to \R$. As above, for $q=\infty$, this space is called the Marcinkiewicz space. 
 
For the basic properties of Lorentz spaces we refer to~\cite[Chapter 4, Section 4]{BeSa88}. Note that spaces $\Lambda^{p, q}$ and $\Lambda^{[p, q]}$ are quasi-normed, rearrangement invariant Banach function spaces. Space $\Lambda^{p, q}$ is a Banach function space if and only if $1 \leq q \leq p < \infty$, while $\Lambda^{[p, q]}$ is a Banach function space if and only if $1 \leq p < \infty$ and $1 \leq q \leq \infty$.  Observe that $\Lambda^{p, p} = L^p$, i.e., standard Lebesgue spaces are Lorentz spaces. 

It holds that $\Lambda^{[p, q]} \subseteq \Lambda^{p, q}$ for arbitrary $p, q$. By~\cite[Chapter 4, Lemma 4.5]{BeSa88}, the quasi-norms $||\cdot||_{\Lambda^{[p, q]}}$ and $||\cdot||_{\Lambda^{p, q}}$ are equivalent for $p > 1$, $q \geq 1$, which means that \begin{equation}\label{eq:loreq}
    \Lambda^{[p, q]} = \Lambda^{p, q} \qquad \text{for $p > 1$, $q \geq 1$}\,.
\end{equation} Lorentz spaces also satisfy the embedding $\Lambda^{p, q_2} \subseteq \Lambda^{p, q_1}$ for $q_1 \geq q_2$. In fact, by~\cite[Chapter 4, Proposition 4.2]{BeSa88} there exists a constant $C = C(q_1, q_2)$ such that
\begin{equation}\label{eq:lorinq}
  ||f||_{\Lambda^{p, q_1}} \leq C||f||_{\Lambda^{p, q_2}}\,.
\end{equation}  
Moreover, due to \cite[page 217]{BeSa88}, we know that if $\Omega \subset \rn$ is such that $|\Omega| < \infty$, then
\begin{equation}\label{eq:loremb}
    \text{if $0 < p_1 < p_2 \leq \infty$ and $0 < q_1, q_2 \leq \infty$, then } \Lambda^{p_2, q_2}(\Omega) \subseteq \Lambda^{p_1, q_1}(\Omega)\,.
\end{equation}

\section{Properties of potential $\Waps$. Proofs of Theorem~\ref{theo:boundedness-psi}, \ref{theo:equivalence}, and~\ref{theo:equivalence-operators}, and direct consequences.}\label{sec:main-proofs}

Let us recall that potential $\Waps$ is given by \eqref{wolff-potential-psi}.  We denote the truncated potential by
\begin{equation}\label{eq:WR-psi} \WapsR f(x) := \int_{0}^{R}r^{\alpha - 1} \psi\left(r^{\alpha-n}{\int_{B(x,r)}|f(y)|\,dy}\right)dr\,. \end{equation}
Before proofs of our main results on one-dimensional reduction, let us present basic properties of $\WapsR$.
\begin{lem}\label{lem:WR-psi-M}Let $n\in\N$ and $\alpha\in (0,n)$.  Assume $\psi:\rp\to\rp$ is non-decreasing, then for every $x\in\rn$ it holds that
\begin{equation*}
\WapsR f(x) \leq \tfrac{1}{\alpha} R^{\alpha}\psi\left(\omega_n^{1 - \frac{\alpha}{n}} M_{\alpha}f(x)\right)\qquad\text{for  $f\in L^1_{loc}(\rn)\,$}.
\end{equation*} 
\end{lem}
\begin{proof}
We observe that by the very definition of the maximal function~\eqref{eq:maxop} we have
\begin{align*}
\WapsR f (x) &= \int_{0}^{R}r^{\alpha - 1}\psi\left(r^{\alpha-n}\int_{B(x,r)}|f(y)|dy\right)dr\\
&\leq \int_{0}^{R}r^{\alpha - 1}\psi\left(\omega_n^{1 - \frac{\alpha}{n}} M_{\alpha}f(x)\right)dr =\frac{R^{\alpha}}{\alpha} \psi\left(\omega_n^{1 - \frac{\alpha}{n}} M_{\alpha}f(x)\right)\,.
\end{align*}
\end{proof}

\begin{lem}\label{lem:WR-psi}Let $n\in\N$ and $\alpha\in (0,n)$.  Suppose $\psi:\rp\to\rp$ is a non-decreasing, left-continuous function and  $f\in L^1_{loc}(\rn)$, then there exists $c=c(\alpha,n)>0$ such that for every $t>0$ it holds that
\begin{equation}
(\WapsR f)^*(t) \leq \tfrac{1}{\alpha} R^{\alpha}\psi\left(c\sup_{s > t}  s^{\frac{\alpha}{n}}f^{**}(s)\right)\qquad\text{for  $f\in L^1_{loc}(\rn)\,$}.
\end{equation} 
\end{lem}
\begin{proof}
Function $t \mapsto \tfrac{1}{\alpha}R^\alpha \psi(c\,t)$ is a non-decreasing, left-continuous function for a fixed $c> 0$. Therefore, Lemma~\ref{lem:psi-rearr} implies that
\[ \left(\frac{R^{\alpha}}{\alpha}\psi\left(\omega_n^{1 - \frac{\alpha}{n}} M_{\alpha}f\right)\right)^{*}(t) = \frac{R^{\alpha}}{\alpha}\psi\left(\omega_n^{1 - \frac{\alpha}{n}}\left( M_{\alpha}f\right)^{*}(t)\right). \]
By using~Lemma~\ref{lem:WR-psi-M} in conjunction with~\eqref{remax} we obtain
\[ (\WapsR f)^{*}(t) \leq \tfrac{1}{\alpha}R^{\alpha}\psi\left(\omega_n^{1 -\frac{\alpha}{n}} \left(M_{\alpha}f\right)^{*}(t)\right)\leq  \tfrac{1}{\alpha}R^{\alpha}\psi\left(C_M\omega_n^{1 - \frac{\alpha}{n}} \sup_{s > t} s^{\frac{\alpha}{n}}f^{**}(s)\right). \]
\end{proof}

Now we shall investigate the operator $f \mapsto \left(\Waps f - \WapsR f\right)$. We use the following fact similar to~\cite[Lemma~2.5]{CYZG-Wolff}.
\begin{lem}\label{lem:Wrest-psi} Let $n\in\N$, $\alpha\in [0,n]$ and let $\psi:\rp\to\rp$ be a non-decreasing function. There exist constants $c_1 = c_1(\alpha, n) > 0$, $c_2 = c_2(\alpha, n) > 0$ such that for every $R > 0$  the following inequality holds true
\begin{equation}\label{eq:Wrest-psi}
    \Waps f(x) - \WapsR f(x)  \leq c_1 \int_{\omega_nR^n}^{\infty} s^{\frac{\alpha}{n}-1}\psi\left(c_2 s^{\frac{\alpha}{n}}f^{**}(s)\right)ds\qquad\text{for  $f\in L^1_{loc}(\rn)\,$}.
\end{equation}
\end{lem}
\begin{proof}  Observe that by~\eqref{eq:decinq} for any $x \in \mathbb{R}^n$ and $r > 0$ it holds that
\[ \int_{B(x, r)} |f(y)|\,dy \leq \int_{0}^{|B(x, r)|} f^{*}(s)\,ds = \int_{0}^{\omega_nr^n} f^{*}(s)\,ds\,. \]
Therefore
\begin{align*}
     \Waps f(x) - \WapsR f(x)&= \int_{R}^{\infty}r^{\alpha - 1} \psi\left(r^{\alpha-n}\int_{B(x, r)}|f(y)|\,dy\right)\,dr \leq \int_{R}^{\infty}r^{\alpha - 1}\psi\left(r^{\alpha-n}\int_{0}^{\omega_n r^n} f^{*}(t)\,dt\right)dr \\
& = \int_{R}^{\infty}r^{\alpha - 1} \psi\left(\omega_nr^{\alpha}f^{**}(\omega_nr^n)\right)\,dr = \frac{\omega_n^{{-\frac{\alpha}{n}}}}{n} \int_{\omega_n R^n}^{\infty} s^{\frac{\alpha}{n}-1}\psi\left(\omega_n^{1-\frac{\alpha}{n}}s^{\frac{\alpha}{n}}f^{**}(s)\right)\,ds\,, 
\end{align*}
where in the last equality we used a substitution $s = \omega_n r^n$.
\end{proof}

We are in position to prove Theorem~\ref{theo:boundedness-psi}.

\begin{proof}[Proof of Theorem~\ref{theo:boundedness-psi}] 
Observe that the right-hand side of~\eqref{eq:Wrest-psi} does not depend on $x$. We fix arbitrary $R>0$. Combining Lemmas~\ref{lem:WR-psi} and~\ref{lem:Wrest-psi} gives us the following estimate
\begin{equation} 
    (\Waps f)^{*}(t) \leq \tfrac{1}{\alpha}R^{\alpha}\psi\left(c\sup_{s > t} s^{\frac{\alpha}{n}}f^{**}(s)\right) + c_1\int_{\omega_nR^n}^{\infty} s^{\frac{\alpha}{n}-1}\psi(c_2s^{\frac{\alpha}{n}}f^{**}(s))ds\,.
\end{equation} 
By substituting $t = \omega_nR^n$ we obtain that
\begin{equation}\label{eq:Wt}
    (\Waps f)^{*}(t) \leq  \frac{1}{\alpha \omega_n^{\frac{\alpha}{n}}}t^{\frac{\alpha}{n}}\psi\left(c\sup_{s > t} s^{\frac{\alpha}{n}}f^{**}(s)\right) + c_1\int_{t}^{\infty} s^{\frac{\alpha}{n}-1}\psi(c_2s^{\frac{\alpha}{n}}f^{**}(s))ds\,.
\end{equation} 
To complete the proof of~\eqref{eq:West-psi}, we shall show that the first summand of~\eqref{eq:Wt} is always lesser (up to a constant) than the second one. In fact, it suffices to show that
\[ t^{\frac{\alpha}{n}}\psi\left(c\tau^{\frac{\alpha}{n}}f^{**}(\tau)\right) \leq  2^{1-\frac{\alpha}{n}}\int_{t}^{\infty} s^{\frac{\alpha}{n}-1}\psi\left(2cs^{\frac{\alpha}{n}}f^{**}(s)\right)\,ds \]
for every $\tau > t$. Observe that due to~\eqref{max-rear-prop} we have
\begin{align*}
    t^{\frac{\alpha}{n}}\psi(c\tau^{\frac{\alpha}{n}}f^{**}(\tau))&\leq \tau^{\frac{\alpha}{n}}\psi(c\tau^{\frac{\alpha}{n}}f^{**}(\tau)) =  \int_{\tau}^{2\tau} \tau^{\frac{\alpha}{n}-1}\psi(c\tau^{\frac{\alpha}{n}}f^{**}(\tau))\,ds\\
    &\leq \int_{\tau}^{2\tau} 2^{1-\frac{\alpha}{n}}s^{\frac{\alpha}{n}-1}\psi(cs^{\frac{\alpha}{n}}f^{**}(\tfrac12s))\,ds \leq 2^{1-\frac{\alpha}{n}}\int_{t}^{\infty} s^{\frac{\alpha}{n}-1}\psi(2cs^{\frac{\alpha}{n}}f^{**}(s))\,ds\,.
\end{align*}
This ends the proof of the upper bound.

\medskip

To prove sharpness, we take any nonnegative radially decreasing function $f$. We notice that
\begin{equation*} \Waps f(x) = \int_{0}^{\infty}r^{\alpha - 1} \psi\left(r^{\alpha-n}\int_{B(x, r)}f(y)\,dy\right)dr \geq \int_{2|x|}^{\infty}r^{\alpha - 1} \psi\left(r^{\alpha-n}\int_{B(x, r)}f(y)\,dy\right)dr\,. \end{equation*}
For $r > |x|$ it holds that $B(0, r - |x|) \subseteq B(x, r)$. In turn, we have
\begin{equation}\label{eq:lipiec1} \Waps f(x) \geq \int_{2|x|}^{\infty}r^{\alpha - 1} \psi\left(r^{\alpha-n}\int_{B(0, r - |x|)}f(y)\,dy\right)dr = \int_{2|x|}^{\infty}r^{\alpha - 1} \psi\left(r^{\alpha-n}\int_{0}^{\omega_n(r-|x|)^n}f^*(s)\,ds\right)dr\,, \end{equation}

where we used~\eqref{eq:deceq} in the last equality. As for $r > 2|x|$ it holds that $ \frac{1}{2}r \leq r - |x| \leq r$, we have
\begin{equation}\label{eq:lipiec2} r^{\alpha-n}\int_{0}^{\omega_n(r-|x|)^n}f^*(s)\,ds \geq 2^{{\alpha}-n}\omega_n(r-|x|)^{\alpha}f^{**}(\omega_n(r-|x|)^n)\,. \end{equation}
By substituting~\eqref{eq:lipiec2} into~\eqref{eq:lipiec1} we obtain inequality 
\begin{align*}
    \Waps f(x) &\geq \int_{2|x|}^{\infty} r^{\alpha - 1}\psi \left( 2^{{\alpha}-n}\omega_n(r-|x|)^{\alpha}f^{**}(\omega_n(r-|x|)^n) \right)dr\\
    &\geq C_{\bW, 3}\int_{\omega_n|x|^n}^{\infty} s^{\frac{\alpha}{n}-1}\psi\left(C_{\bW, 4}s^{\frac{\alpha}{n}}f^{**}(s)\right)\,ds\,,
\end{align*}
where the last inequality comes from substitution $s = \omega_n (r-|x|)^n$ and holds for some constants $C_{\bW, 3} = C_{\bW, 3}(\alpha, n) > 0$ and $C_{\bW, 4} = C_{\bW, 4}(\alpha, n) > 0$.
The right-hand side of the last display is nonnegative and radially decreasing function of $x$, so it satisfies~\eqref{eq:radec}. Using~\eqref{eq:monrea}, we obtain
\begin{equation*} (\Waps f)^*(t) \geq C_{\bW, 3}\int_{t}^{\infty} s^{\frac{\alpha}{n}-1}\psi(C_{\bW, 4}s^{\frac{\alpha}{n}}f^{**}(s))\,ds\,. \end{equation*}
This gives us~\eqref{eq:Wsharp}.
\end{proof}

We make use of Theorem~\ref{theo:boundedness-psi} to prove boundedness criteria.
\begin{proof}[Proof of Theorem~\ref{theo:equivalence}]
We start with the implication \textit{(i)} $ \Rightarrow $ \textit{(ii)}. Having any $f \in X(\Omega)$, we take \mbox{$f^*\in \overline{X}(0,|\Omega|)$}.  Then Theorem~\ref{theo:boundedness-psi} \blue{gives} us that
\begin{align*}
 ||\Waps f||_{Y(\Omega)} &\leq {
 ||(\Waps f{)}^*||_{\overline{Y}(0,|\Omega|)} } \leq  \left| \left|C_{\bW,1} \int_{{t}}^{\infty} s^{\frac{\alpha}{n}-1}\psi\left(C_{\bW,2} s^{\frac{\alpha}{n}-1}\int_{0}^{s}f^*(y)dy\right)ds \right| \right|_{\overline{Y}(0, |\Omega|)}\,.
\end{align*} 
We can estimate the right-hand side above using  \textit{(i)} to get
\begin{align*}
 ||\Waps f||_{Y(\Omega)} \leq cC_{\bW, 1} ||C_{\bW, 2}f^*||_{\overline{X}(0, |\Omega|)} = cC_{\bW, 1}C_{\bW, 2} ||f||_{X(\Omega)}\,. 
\end{align*}
Hence, the implication is proved. \newline

Conversely, in order to prove  \textit{(ii)} $ \Rightarrow $ \textit{(i)} we take any nonnegative $\phi \in \overline{X}(0, \infty)$.   Let $f(x) = \tfrac{1}{C_{\bW, 4}} \phi^*(\omega_n|x|^n)$, so that $f$ is a nonnegative radially decreasing function and $f^*(t) = \tfrac{1}{C_{\bW, 4}}\phi^*(t)$. Notice that for every $r > 0$ it holds that
\[ \int_{0}^{r} \phi(s)ds \leq \int_{0}^{r} \phi^{*}(s)ds\,, \]
which implies that
\begin{align*}
& \left| \left| \int_{{t}}^{\infty} s^{\frac{\alpha}{n}-1}\psi\left(s^{\frac{\alpha}{n}-1}\int_{0}^{s}\phi(y)dy\right)ds \right| \right|_{\overline{Y}(0, \infty)}
\leq \left| \left|  \int_{{t}}^{\infty} s^{\frac{\alpha}{n}-1}\psi\left(s^{\frac{\alpha}{n}-1}\int_{0}^{s}\phi^*(y)dy\right)ds \right| \right|_{\overline{Y}(0, \infty)}=:I_1\,.
\end{align*}
Since $f$ is radially decreasing and $\Omega=\rn$, Theorem~\ref{theo:boundedness-psi} implies that
\begin{equation*} C_{\bW, 3} \int_{{t}}^{\infty} s^{\frac{\alpha}{n}-1}\psi\left(s^{\frac{\alpha}{n}}C_{\bW, 4}f^{**}(s)\right)ds\leq (\Waps f)^{*}(t)\,. \end{equation*}
Hence, we infer that
\begin{align*}
I_1 \leq \tfrac{1}{C_{\bW, 3}}||(\Waps f)^*||_{\overline{Y}(0, \infty)} = \tfrac{1}{C_{\bW, 3}}||\Waps f||_{Y(\rn)}\,.\end{align*}
Taking into account \textit{(ii)} and collecting the above observations, we see that
 \begin{align*}
     \left| \left|\int_{{t}}^{\infty} s^{\frac{\alpha}{n}-1}\psi\left(s^{\frac{\alpha}{n}-1}\int_{0}^{s}\phi(y)dy\right)ds \right| \right|_{\overline{Y}(0, \infty)}&\leq \tfrac{1}{C_{\bW, 3}}||\Waps f||_{Y(\rn)}  \leq \tfrac{1}{C_{\bW, 3}}c ||f||_{X(\rn)}\\ &=\tfrac{1}{C_{\bW, 3}C_{\bW, 4}}c||\phi^*||_{\overline{X}(0, \infty)}=\tfrac{1}{C_{\bW, 3}C_{\bW, 4}} c ||\phi||_{\overline{X}(0, \infty)}\,,
 \end{align*}
which was to be proven.
\end{proof}

\begin{proof}[Proof of Theorem~\ref{theo:equivalence-operators}]
We start with the implication \textit{(i)} $ \Rightarrow $ \textit{(ii)}. Let $k>0$ be arbitrarily fixed. Having any $f \in X(\Omega)$, we take $f^*\in \overline{X}(0,|\Omega|)$. Then~\eqref{oT-def}, \eqref{oTY-monot}, and Theorem~\ref{theo:boundedness-psi} give us that
\begin{align*}
 \blue{\mathcal{T}_Y(k\Waps f)} &\blue{\leq \oT_{\overline{Y}}(k(\Waps f)^*) }\leq \oT_{\overline{Y}}\left(kC_{\bW,1} \int_{{t}}^{\infty} s^{\frac{\alpha}{n}-1}\psi\left(C_{\bW,2}s^{\frac{\alpha}{n}-1}\int_{0}^{s} f^*(y)\,dy\right)ds\right) \,.
\end{align*}
By assumption \textit{(i)}, for some $c > 0$ we have
\begin{align*}
 \mathcal{T}_Y(k\Waps f){\leq} \oT_{\overline{X}}(c C_{\bW, 2}f^*) =  \mathcal{T}_{{X}}(c C_{\bW, 2}f) \,.
\end{align*}
Hence, the implication is proven. \newline

Conversely, in order to prove  \textit{(ii)} $ \Rightarrow $ \textit{(i)} we take any nonnegative $\phi \in \overline{X}(0, \infty)$ and $k>0$.  Let $f(x) =\tfrac{1}{C_{\bW, 4}} \phi^*(\omega_n|x|^n)$, so that $f$ is a nonnegative radially decreasing function and $f^*(t) =\tfrac{1}{C_{\bW, 4}} \phi^*(t)$. Notice that for every $r > 0$ it holds that
\[ \int_{0}^{r} \phi(s)ds \leq \int_{0}^{r} \phi^{*}(s)ds\,, \]
which in the view of~\eqref{oTY-monot} implies that
\[ \oT_{\overline{Y}}\left( k\int_{{t}}^{\infty} s^{\frac{\alpha}{n}-1}\psi\left(s^{\frac{\alpha}{n}-1}\int_{0}^{s}\phi(y)dy\right)ds\right)  \leq  \oT_{\overline{Y}}\left(k \int_{{t}}^{\infty} s^{\frac{\alpha}{n}-1}\psi\left(s^{\frac{\alpha}{n}-1}\int_{0}^{s}\phi^*(y)dy\right)ds\right) =: I_2\,.\]
Since $f$ is radially decreasing and $\Omega=\rn$, Theorem~\ref{theo:boundedness-psi} implies that
\[  C_{\bW, 3}\int_{{t}}^{\infty} s^{\frac{\alpha}{n}-1}\psi\left(s^{\frac{\alpha}{n}}C_{\bW, 4}f^{**}(s)\right)ds\leq (\Waps f)^{*}(t)\,. \]
Hence, by~\eqref{oTY-monot} and~\eqref{oT-def} we infer  that
\begin{align*}
I_2 &\leq  \oT_{\overline{Y}} (\tfrac{1}{C_{\bW, 3}}k (\Waps f)^* )=   \mathcal{T}_{{Y}}(\tfrac{1}{C_{\bW, 3}}k \Waps f)\,.\end{align*}
Taking into account \textit{(ii)} and collecting the above observations we see that for some constant $c > 0$ it holds that
 \begin{align*}
      \oT_{\overline{Y}}\left( k\int_{{t}}^{\infty} s^{\frac{\alpha}{n}-1}\psi\left(s^{\frac{\alpha}{n}-1}\int_{0}^{s}\phi(y)dy\right)ds\right) &\leq  \mathcal{T}_{{X}}\big(cf\big) = \oT_{\overline{X}}(cC_{\bW, 4}\phi^*) =  \oT_{\overline{X}}(cC_{\bW, 4}\phi)\,,
 \end{align*}
which was to be proven.
\end{proof}

We have the following consequences of Theorem~\ref{theo:boundedness-psi}.
\begin{coro}
\label{coro:full-Wolff-bdd} Suppose $\psi:\rp\to\rp$ is a non-decreasing function, $n\geq 1$, and $\alpha\in(0,n)$. Then
\begin{itemize}
    \item [(i)] if $\sigma>1$, $\beta,\vr>0$, $\alpha < \frac{\vr n}{\sigma \beta + \sigma \vr}$, and $\psi^\beta(t)\leq ct^\vr$, then for $\gamma=\frac{\beta\sigma n}{\vr n-\alpha\sigma\beta - \alpha \sigma \vr}$ it holds \[\Waps :\Lambda^{\sigma,\vr}(\rn)\to\Lambda^{\gamma,\beta}(\rn)\,;\]
    \item[(ii)]  if $\vr > \frac{\alpha}{n-\alpha}$ and $\psi(t)\leq ct^\vr$, then it holds \[\Waps :L^1(\rn)\to\Lambda^{\frac{n}{n\vr-\alpha(\vr+1)},\infty}(\rn)\,.\]
\end{itemize}
\end{coro}
Note that this result in the case of $\psi(t)=t^{1/(p-1)}$ yields the result of \cite[Theorem~3.1, (i) and (ii)]{Ci-pot} with $\Wap$ in the place of $\cV_{\alpha,p}$.
\begin{proof}
{\it (i)} Let us notice that by Theorem~\ref{theo:boundedness-psi}, Hardy inequality, and the fact that $\psi^\beta(t)\leq ct^\vr$, we have
\begin{align*}
    \int_0^\infty\left(t^\frac{1}{\gamma}(\Waps f)^*(t)\right)^\beta\frac{dt}{t}&\leq 
    C\int_0^\infty\left(t^\frac{1}{\gamma}\int_t^\infty s^{\frac{\alpha}{n}-1}\psi(c s^\frac{\alpha}{n} f^{**}(s))\,ds\right)^\beta\frac{dt}{t} =C
    \int_0^\infty t^{\frac{\beta}{\gamma}-1}\left(\int_t^\infty s^{\frac{\alpha}{n}-1}\psi(c s^\frac{\alpha}{n} f^{**}(s))\,ds\right)^{\beta} \,{dt}\\
    &\leq C\int_0^\infty t^{{\beta}\left(\frac{1}{\gamma}+\frac{\alpha}{n}\right)} \psi^\beta(c t^\frac{\alpha}{n} f^{**}(t)) \frac{dt}{t}\leq Cc \int_0^\infty t^{{\beta}\left(\frac{1}{\gamma}+\frac{\alpha}{n}\right)}  t^\frac{\alpha\vr}{n} (f^{**}(t))^\vr \frac{dt}{t}=C\|f\|^\vr_{\Lambda^{[\sigma,\vr]}}\,,
\end{align*}
where the last equality comes from the definition of $\gamma$. Note that since $\sigma>1$, we have $\Lambda^{\sigma,\vr}=\Lambda^{[\sigma,\vr]}.$ \newline

{\it (ii)} Theorem~\ref{theo:boundedness-psi} and $\psi(t)\leq ct^\vr$ imply that
\[(\Waps f)^*(t)\leq C\int_{t^\frac{\alpha}{n}}^\infty \psi(c r f^{**}(r^\frac{n}{\alpha}))\,dr\leq C\int_{t^\frac{\alpha}{n}}^\infty \psi(c r^{1-\frac{n}{\alpha}} \|f\|_{L^1})\,dr \leq Cc \int_{t^\frac{\alpha}{n}}^\infty  r^{\vr(1-\frac{n}{\alpha})} \|f\|_{L^1}^\vr\,dr\leq C\|f\|^\vr_{L^1}  t^{\frac{\alpha -n}{n}\vr+\frac{\alpha}{n}}\,.\]
Therefore\begin{align*}
   \sup_{t>0} t^\frac{n\vr-\alpha(\vr+1)}{n}(\Waps f)^*(t)& \leq c\|f\|_{L^1}^\vr\,.
\end{align*}
\end{proof}

\begin{lem} Suppose $\psi:\rp\to\rp$ is a non-decreasing function, $n\geq 1$, and $\alpha\in [0,n]$, then there exists a~constant $c=c(n,\alpha)>0$ such that
 \label{lem:Wolff:Lambda-Linfty}\[
 \|\WapsR f\|_{L^\infty(\R^n)}\leq c \int_0^{\omega_n R^n}  t^{\frac{\alpha}{n}-1} \psi {\left(\omega_n^{1-\frac{\alpha}{n}}t^\frac{\alpha}{n}f^{**}(t)\right)}\,{dt}\,.\]
Moreover, $\Waps$ maps a Lorentz-type space of functions such that $ \int_0^{\infty}  t^{\frac{\alpha}{n}-1} \psi {\left(\omega_n^{1-\frac{\alpha}{n}}t^\frac{\alpha}{n}f^{**}(t)\right)}\,{dt}<\infty$ into $L^\infty(\rn)$ and \[ 
 \|\Waps f\|_{L^\infty(\R^n)}\leq c \int_0^{\infty}  t^{\frac{\alpha}{n}-1} \psi {\left(\omega_n^{1-\frac{\alpha}{n}}t^\frac{\alpha}{n}f^{**}(t)\right)}\,{dt}\,.\] 
\end{lem}

\begin{proof}
By~\eqref{eq:decinq} we have
\begin{align*}
    ||\WapsR{f}||_{L^{\infty}(\R^n)} \leq \int_{0}^{R} r^{\alpha - 1} \psi \left( r^{\alpha - n}\omega_nr^nf^{**}(\omega_nr^n) \right)\,dr &= \int_{0}^{R} r^{\alpha - 1}\psi \left( \omega_n r^{\alpha}f^{**}(\omega_n r^n) \right)\,dr\\
    &\leq c\int_{0}^{\omega_n R^n} t^{\frac{\alpha}{n} - 1}\psi(\omega_n^{1 - \frac{\alpha}{n}}t^{\frac{\alpha}{n}}f^{**}(t))\,dt\,.
\end{align*}
\end{proof}

\section{Potentials given by a doubling function and regularity to quasilinear PDEs}\label{sec:appl}
The aim of this section is to present some applications of Theorems~\ref{theo:boundedness-psi} and~\ref{theo:equivalence} to the theory of regularity. In fact, knowing the estimates for solutions to problems with doubling Orlicz growth expressed with use of generalized Wolff potential, we can precisely transfer regularity from data to solutions. The regularity we infer was not known in that delicate scales. We start with some more definitions, then we pass to boundedness of $\WaG$ and regularity.
\subsection{$N$-functions}\label{ssec:N-functions}
Basic references for what follows are~\cite{adamsfournier, rao-ren}. In this section let  $G\in C(\rp)\cap C^1((0,\infty))$ and assume that \begin{itemize}
    \item $G$ is an $N$-function -- it is a convex, continuous,  and such that $G(0)=0$, $\lim_{t \to 0}{G(t)}/{t}=0$ and $\lim_{t \to \infty}{G(t)}/{t}=\infty$;
    \item $g:\rp\to\rp$ is a non-decreasing function such that $G(t)=\int_0^t g(s)\,ds\,$;
    \item $G \in \Delta_2 \cap \nabla_2$, which is equivalent to the fact that there exist constants $i_G,s_G>1$ satisfying \begin{equation}\label{iG-sG}
1<  i_G=\inf_{t>0}\frac{tg(t)}{G(t)}\leq \sup_{t>0}\frac{tg(t)}{G(t)}=s_G<\infty\,.\end{equation}
\end{itemize}

The Young conjugate $\wt{G}$ to an $N$-function $G:\rp\to\rp$  is given by $\wt{G}(s):=\sup_{t>0}(s t-G(t)).$ Then $\wt G$ is an $N$-function too. 

We say that a function $G:\rp\to\rp$ satisfies $\Delta_2$-condition (denoted $G\in\Delta_2$), if there exists $c_{\Delta_2}>0$ such that $G(2t)\leq c_{\Delta_2}G(t)$ for $t>0$. We say that $G$ satisfy $\nabla_2$-condition if $\wt{G}\in\Delta_2.$ Note that it is possible that $G$ satisfies only one of the conditions $\Delta_2/\nabla_2$. For instance, when $\blue{G(t) = (1+t)\log(1+t)-t}$, its complementary function is  $\widetilde{G}(s)= \blue{\exp(s)-s-1}$. Then $G\in\Delta_2$ and $\widetilde{G}\not\in\Delta_2$.  Note that $G\not\in\Delta_2$ can be trapped between two power-type functions with arbitrarily close powers, see~\cite{CGZG,BDMS}. If $G\in C([0, \infty))\cap C^1((0,\infty))$ is an $N$-function and $G' = g$, then $G\in\Delta_2\cap\nabla_2$ holds if and only if~\eqref{iG-sG}. Moreover, then
\begin{equation}\label{comp-i_G-s_G} 
\frac{G(t)}{t^{i_G}}\quad\text{is non-decreasing}\qquad\text{and}\qquad\frac{G(t)}{t^{s_G}}\quad\text{is non-increasing}\,,
\end{equation} 
and $i_GG(t)\leq g(t)t\leq s_GG(t).$ Note that $g(t)t$ is comparable to $G(t)$, so function $t \mapsto \frac{g(t)}{t^{i_G - 1}}$ is comparable to a~non-decreasing function and $t \mapsto \frac{g(t)}{t^{s_G - 1}}$ is comparable to a non-increasing function. Observe that $i_{\wt{G}} = \tfrac{s_G}{s_G - 1}$ and $s_{\wt{G}} = \tfrac{i_G}{i_G - 1}$.

\subsection{Properties of $\WaG$}\label{sec:WaG}

We start with a consequence of Theorem~\ref{theo:boundedness-psi}. It follows directly from this theorem and doubling properties of~$g^{-1}$. In the further parts of this section, we make use of this result to establish various results on the boundedness of $\WaG$.

\begin{prop} \label{prop:boundedness-G}
Let $n\geq 1$, $\alpha\in(0,n)$, and let $G\in C(\rp)\cap C^1((0,\infty))$ be an $N$-function such that $G\in\Delta_2\cap\nabla_2$. Then there exists $C_{\cW}=C_\cW(\alpha,i_G, s_G, n) > 0$ such that
\begin{equation}\label{eq:West}(\WaG f)^{*}(t) \leq C_\cW \int_{{t}}^{\infty} s^{\frac{\alpha}{n}-1}g^{-1}\left(s^{\frac{\alpha}{n}}f^{**}(s)\right)ds\end{equation} whenever $f\in L^1_{loc}(\rn)$ satisfies~\eqref{skoncz-wart}.\\
The result
is sharp, in the sense that there exists a positive constant $c_\cW =c_\cW (\alpha, i_G, s_G, n) > 0$ such that, if $f\in L^1_{loc}(\rn)$ is nonnegative and radially decreasing, then
\begin{equation}\label{eq:Wsharp}
c_\cW\int_{{t}}^{\infty} s^{\frac{\alpha}{n}-1}g^{-1}\left(s^{\frac{\alpha}{n}}f^{**}(s)\right)ds \leq  (\WaG f)^*(t).
\end{equation}
\end{prop}

We use Proposition~\ref{prop:boundedness-G} to prove the following estimates.

\begin{prop}\label{prop:WaGest} Suppose  $G\in C(\rp)\cap C^1((0,\infty))$ is an $N$-function such that $G\in\Delta_2\cap\nabla_2$. The following estimates hold true.
\begin{enumerate}
    \item [(i)] For any $\alpha \in (0, \frac{n}{s_G})$ and $q \in (0, 1]$ there exists a constant $C = C(\alpha, n, i_G, s_G, q) > 0$ such that
    \begin{equation}\label{eq:WaG1}
        t^{1 - \frac{\alpha}{n}}g(t^{-\frac{\alpha}{n}} (\WaG{f})^{*}(t) ) \leq C||f||_{\Lambda^{1, q}(\R^n)}\,
    \end{equation}
    for every $f \in \Lambda^{1, q}(\rn)$.
    \item [(ii)] For any $\beta > 1$, $\alpha \in (0, \frac{n}{\beta s_G})$ and $q \in (0, \infty]$ there exists a constant $C = C(\beta, \alpha, n, i_G, s_G, q) > 0$ such that
    \begin{equation}\label{eq:WaG2}
        t^{\frac{1}{\beta} - \frac{\alpha}{n}}g(t^{-\frac{\alpha}{n}} (\WaG{f})^{*}(t) ) \leq C||f||_{\Lambda^{\beta, q}(\R^n)}\,
    \end{equation}
    
    for every $f \in \Lambda^{\beta, q}(\rn)$.
    
\end{enumerate}
\end{prop}

\begin{proof}
We start with the proof of~\eqref{eq:WaG1} for $q=1$. Combining Proposition~\ref{prop:boundedness-G} and Lemma~\ref{lem:gpow} give us that
\begin{align*}
    \left( \WaG f \right)^{*}(t) &\leq C_{\cW} \int_{t}^{\infty} s^{\frac{\alpha}{n} - 1}g^{-1} \left( s^{\frac{\alpha}{n} - 1} \int_{0}^{s} f^{*}(r)dr \right)ds
    \leq C_{\cW} \int_{t}^{\infty} s^{\frac{\alpha}{n} - 1}g^{-1} \left( s^{\frac{\alpha}{n} - 1} ||f||_{L^1(\R^n)} \right)ds \\
    &= \tfrac{n}{\alpha} C_{\cW} \int_{t^{\frac{\alpha}{n}}}^{\infty} g^{-1} \left(s^{1 - \frac{n}{\alpha}}  ||f||_{L^1(\R^n)}  \right)ds
    \leq \tfrac{n}{\alpha} C_{\cW}c_{G, 2} t^{\frac{\alpha}{n}}g^{-1}\left( t^{\frac{\alpha}{n} - 1}||f||_{L^1(\R^n)} \right),
\end{align*}
    which is~\eqref{eq:WaG1} with $q = 1$. Using inequality~\eqref{eq:lorinq} gives us the desired result. \newline
    
    The proof of ~\eqref{eq:WaG2} is analogous. We consider $q=\infty$. Let us notice that due to the fact that $\beta>1$ and~\eqref{eq:loreq} we have
\begin{align*}
    \left( \WaG f \right)^{*}(t) &\leq C_{\cW} \int_{t}^{\infty} s^{\frac{\alpha}{n} - 1}g^{-1} \left( s^{\frac{\alpha}{n}} f^{**}(s) \right)ds
    \leq C_{\cW} \int_{t}^{\infty} s^{\frac{\alpha}{n} - 1}g^{-1} \left( s^{\frac{\alpha}{n} - \frac{1}{\beta}} ||f||_{\Lambda^{\beta, \infty}(\rn)} \right)ds \\
    &= \tfrac{n}{\alpha} C_{\cW} \int_{t^{\frac{\alpha}{n}}}^{\infty} g^{-1} \left(s^{1 - \frac{n}{\alpha \beta}}  ||f||_{\Lambda^{\beta, \infty}(\rn)}  \right)ds
    \leq \tfrac{n}{\alpha} C_{\cW}c_{G, 2} t^{\frac{\alpha}{n}}g^{-1}\left( t^{\frac{\alpha}{n} - \frac{1}{\beta}}||f||_{\Lambda^{\beta, \infty}(\rn)} \right),
\end{align*}
    which by~\eqref{eq:lorinq} ends the proof.
\end{proof}

\begin{rem}
By taking $G(t) = \frac{t^2}{2}$ we obtain $\WaG = \cI_{2\alpha} = \cV_{\alpha, 2}$ and the result of Proposition~\ref{prop:boundedness-G} coincides with the classical estimate for Riesz potential~\eqref{eq:ONeil} due to~\cite{ONeil}. Inequality~\eqref{eq:WaG1} implies that $\cI_{\alpha} : L^1(\rn) \to \Lambda^{\frac{n}{n-\alpha}, \infty}(\rn)$ is bounded.
\end{rem}

To achieve more transparent results, we restrict our attention to functions with support included in a set of finite measure.

\begin{prop}\label{prop:WgOmegaest}
Let $n \geq 1$, $\alpha \in (0, \tfrac{n}{s_G})$, and $G\in C(\rp)\cap C^1((0,\infty))$ be an $N$-function such that $G\in\Delta_2\cap\nabla_2$. Assume further that $\Omega$ is a measurable set such that $|\Omega| < \infty$. Then there exists a constant $C = C(\alpha, n, i_G, s_G) > 0$ such that for all $t \leq |\Omega|$ it holds that
\begin{equation*}
    (\WaG f)^*(t) \leq C\int_{\frac12 t}^{|\Omega|} s^{\frac{\alpha}{n}-1}g^{-1}\left(s^{\frac{\alpha}{n}}f^{**}(s)\right)ds
\end{equation*}
for every measurable \blue{$f : \Omega \to \mathbb{R}$}.
\end{prop}
Note that it is possible to prove this proposition with the integral over $[ct,|\Omega|)$ for any fixed $c\in(0,1)$.
\begin{proof}By Proposition~\ref{prop:boundedness-G} it suffices to show that
\begin{equation*}
    \int_{|\Omega|}^{\infty} s^{\frac{\alpha}{n} - 1}g^{-1} \left(s^{\frac{\alpha}{n}}f^{**}(s) \right)ds \leq c\int_{\frac{1}{2}t}^{|\Omega|} s^{\frac{\alpha}{n} - 1}g^{-1}\left(s^{\frac{\alpha}{n}}f^{**}(s) \right)\,ds
\end{equation*}
for $t \leq |\Omega|$ and for some constant $c = c(\alpha, n, i_G, s_G) > 0$. By Lemma~\ref{lem:gpow} and since $\alpha\in (0,\frac{n}{s_G})$, it holds that
\begin{align*}
    &\int_{|\Omega|}^{\infty} s^{\frac{\alpha}{n} - 1}g^{-1} \left(s^{\frac{\alpha}{n}}f^{**}(s) \right)ds 
    = \tfrac{n}{\alpha} \int_{|\Omega|^{\frac{\alpha}{n}}}^{\infty} g^{-1} \left(t^{1 - \frac{n}{\alpha}}||f||_{L^1(\Omega)} \right)dt \leq  C|\Omega|^{\frac{\alpha}{n}}g^{-1}\left(|\Omega|^{\frac{\alpha}{n} - 1}||f||_{L^1(\Omega)} \right).
\end{align*}
Observe that since $||f||_{L^1(\Omega)} = |\Omega|f^{**}(|\Omega|)$ 
for $t \leq |\Omega|$ it holds that
\begin{align*}
    |\Omega|^{\frac{\alpha}{n}}g^{-1}\left(|\Omega|^{\frac{\alpha}{n} - 1 }||f||_{L^1(\Omega)} \right) &= 2 \int_{\frac{1}{2}|\Omega|}^{|\Omega|} |\Omega|^{\frac{\alpha}{n} - 1}g^{-1}\left(|\Omega|^{\frac{\alpha}{n}}f^{**}(|\Omega|) \right)\,ds \leq 2\int_{\frac{1}{2}|\Omega|}^{|\Omega|} s^{\frac{\alpha}{n} - 1}g^{-1}\left((2s)^{\frac{\alpha}{n}}f^{**}(s) \right)\,ds\\
    &\leq c\int_{\frac{1}{2}|\Omega|}^{|\Omega|} s^{\frac{\alpha}{n} - 1}g^{-1}\left(s^{\frac{\alpha}{n}}f^{**}(s) \right)\,ds
    \leq c\int_{\frac{1}{2}t}^{|\Omega|} s^{\frac{\alpha}{n} - 1}g^{-1}\left(s^{\frac{\alpha}{n}}f^{**}(s) \right)\,ds\,.
\end{align*}
\end{proof}


We formulate a general result on boundedness of $\WaG$ yielding the reduction principle in the manner of Theorem~\ref{theo:boundedness-psi}.
\begin{theo}\label{theo:equivalence-G} Suppose $n\geq 1$, $G\in C(\rp)\cap C^1((0,\infty))$ is an $N$-function such that $G\in\Delta_2\cap\nabla_2$, $\alpha\in(0,\tfrac{n}{s_G})$, and let $\Omega \subset \rn$ be a measurable \blue{ and open} set such that $|\Omega| < \infty$. Let $X(\Omega)$ and $Y(\blue{\rn})$ be quasi-normed Banach rearrangement invariant spaces and let $h_X, h_Y : \rp \to \rp$ be non-decreasing, left-continuous functions such that $h_X, h_Y \in \Delta_2$. Then the following assertions are equivalent.
\begin{itemize}
\item[$(i)$]There exists a constant $c>0$ such that for every nonnegative function $\phi \in \overline{X}(0,|\Omega|)$ such that
    \begin{equation*} \left\| h_Y \left(\int_{\frac{1}{2}t}^{|\Omega|} s^{\frac{\alpha}{n}-1}g^{-1}\left(s^{\frac{\alpha}{n}-1}\int_{0}^{s}\phi(y)\,dy\right)ds\right) \right\|_{\overline{Y}(0,|\Omega|)} \leq c ||h_X(\phi)||_{\overline{X}(0,|\Omega|)}\,. \end{equation*}
\item[$(ii)$] There exists a constant $c>0$ such that for every $f:\rn\to\R$ with $\supp\,f\subset\Omega$, it holds that
    \[||h_Y \left(\WaG f \right)||_{Y(\blue{\rn})} \leq c ||h_X(|f|)||_{X(\Omega)}\,.\]
\end{itemize}
\end{theo}
\begin{proof}
    The implication \textit{(i)}$\implies$\textit{(ii)} follows from a similar argumentation from the proof of Theorem~\ref{theo:equivalence}, where we use Propositions~\ref{prop:boundedness-G} and~\ref{prop:WgOmegaest} instead of Theorem~\ref{theo:boundedness-psi} and doubling properties of $h_X,h_Y,g^{-1}$ instead of the homogeneity of the norm.

    Let us prove the implication \textit{(ii)}$\implies$\textit{(i)}. We fix  $B(x_0, r) \subseteq \Omega$ for some $x_0 \in \Omega$ and $r > 0$, which is possible since $\Omega$ is open. We take any nonnegative $\phi \in \overline{X}(0,|\Omega|)$.  Let us consider a function $\kappa : [0, |\Omega|) \to [0, \infty]$ such that $\kappa^*(t) = \phi^*(t)\mathds{1}_{(0, |B(x_0, r)|)}(t)$. We notice that for every $t \in [0, |\Omega|)$ it holds that
    \begin{equation*}
        \phi^{**}(t) \leq \phi^{**}\left(\frac{t|B(x_0, r)|}{|\Omega|}\right) = \kappa^{**}\left(\frac{t|B(x_0, r)|}{|\Omega|}\right) \leq \frac{|\Omega|}{|B(x_0, r)|}\kappa^{**}(t)\,, 
    \end{equation*}
    where the first inequality comes from the monotonicity of maximal rearrangement, and the equality from the fact that $\kappa^* = \phi^*$ on $(0, |B(x_0, r)|)$, while the last estimate holds due to the doubling property of the maximal rearrangement~\eqref{max-rear-prop}. Note that $C={|\Omega|}/{|B(x_0, r)|}$ is finite since $|\Omega|<\infty$. The inequality from the last display in conjunction with $h_Y \in \Delta_2$ and $g^{-1}\in\Delta_2$, implies that
    \begin{align}\label{eq:lip1}
        h_Y \left(\int_{\frac{1}{2}t}^{|\Omega|} s^{\frac{\alpha}{n}-1}g^{-1}\left(s^{\frac{\alpha}{n}-1}\int_{0}^{s}\phi(y)\,dy\right)ds\right) &\leq Ch_Y \left(\int_{\frac{1}{2}t}^{|\Omega|} s^{\frac{\alpha}{n}-1}g^{-1}\left(s^{\frac{\alpha}{n}}\kappa^{**}(s)\right)ds\right) \nonumber\\
        &\leq Ch_Y \left(\int_{t}^{\infty} s^{\frac{\alpha}{n}-1}g^{-1}\left(s^{\frac{\alpha}{n}}\kappa^{**}(s)\right)ds\right)\,.
    \end{align}
    Let us now take the function $f$ defined by $f(x) \coloneqq  \kappa^*(\omega_n|x-x_0|^n)$. Note that $f^* = \kappa^*$, $\supp{f} \subseteq B(x_0, r) \subseteq \Omega$, and $f$ is a translation of a radially decreasing function. Therefore, by  Proposition~\ref{prop:boundedness-G}, we have that
    \begin{equation*}
        \int_{t}^{\infty} s^{\frac{\alpha}{n}-1}g^{-1}\left(s^{\frac{\alpha}{n}}\kappa^{**}(s)\right)ds \leq C(\WaG f)^*(t)\,.
    \end{equation*}
    Therefore, by~\eqref{eq:lip1} and assumption $(ii)$, it holds that
    \begin{align*}
        &\left\| h_Y \left(\int_{\frac{1}{2}t}^{|\Omega|} s^{\frac{\alpha}{n}-1}g^{-1}\left(s^{\frac{\alpha}{n}-1}\int_{0}^{s}\phi(y)\,dy\right)ds\right) \right\|_{\overline{Y}(0,|\Omega|)} \leq C\|h_Y\left((\WaG f)^*\right)\|_{\overline{Y}(0,|\Omega|)} \blue{\leq} C\|h_Y(\WaG f)\|_{Y(\blue{\rn})}\\
        &\leq Cc\|h_X(|f|)\|_{X(\Omega)} = Cc\|h_X(f^*)\|_{\overline{X}(0, |\Omega|)} = Cc\|h_X(\kappa^*)\|_{\overline{X}(0, |\Omega|)} \leq Cc\|h_X(\phi^*)\|_{\overline{X}(0, |\Omega|)} = Cc\|h_X(\phi)\|_{\overline{X}(0, |\Omega|)}\,,
    \end{align*}
    which ends the proof.
\end{proof}
%
%
%
%
%

Let us comment on the assumptions of Theorem~\ref{theo:equivalence-G}.
\begin{rem}\rm
    We shall stress that in the proof of Theorem~\ref{theo:equivalence-G} we cannot relax the assumptions $G\in\Delta_2\cap\nabla_2$ as well as extend the range of admissible $\alpha$. The first reason comes from the observation that Theorem~\ref{theo:boundedness-psi} takes the form of  Proposition~\ref{prop:boundedness-G} only when $G\in\Delta_2\cap\nabla_2$ and $\alpha\in(0,\frac{n}{s_G})$. The second reason lies in the proof of Proposition~\ref{prop:WgOmegaest}, where the doubling growth is necessary to restrict the integral on the right-hand side to $[\frac 12 t,|\Omega|)$. Moreover, in \eqref{eq:lip1} we also change the variables and the domain of integration, which requires the use of the fact that $g^{-1}\in\Delta_2$ to keep the final estimate in the needed form.
\end{rem}


Proposition~\ref{prop:WaGest} allows us to prove the following result relating to~\cite[Theorem 3.1]{Ci-pot}.

\begin{prop}\label{prop:gWest}
    Suppose that  $G\in C(\rp)\cap C^1((0,\infty))$ is an $N$-function such that $G\in\Delta_2\cap\nabla_2$ and $\Omega \subset \R^n$ satisfies $|\Omega| < \infty$. Then the following estimates hold true.
    \begin{enumerate} 
        \item [(i)] For any $\alpha \in (0, \frac{n}{s_G})$ and $q \in (0, 1]$ there exists a constant $C = C(\alpha, n, i_G, s_G, q, |\Omega|) > 0$ such that
        \begin{equation}\label{eq:WaGOm1}
            \left| \left| g(\WaG{f}) \right| \right|_{\Lambda^{\frac{n}{n - \alpha i_G},\infty}(\Omega)} \leq C||f||_{\Lambda^{1, q}(\Omega)}\,
        \end{equation}
        for every $f : \Omega \to \R$ such that $f \in \Lambda^{1, q}(\Omega)$;
        \item [(ii)] For any $\beta > 1$, $\alpha \in (0, \frac{n}{\beta s_G})$ and $q \in (0, \infty]$ there exists a constant $C = C(\beta, \alpha, n, i_G, s_G, q, |\Omega|) > 0$ such that
        \begin{equation}\label{eq:WaGOm2}
            \left| \left| g(\WaG{f}) \right| \right|_{\Lambda^{\frac{\beta n}{n - \beta \alpha i_G},\infty}(\Omega)} \leq C||f||_{\Lambda^{\beta, q}(\Omega)}\,
        \end{equation}
        for every $f : \Omega \to \R$ such that $f \in \Lambda^{\beta, q}(\Omega)$.
        \item [(iii)] For any $\beta > 1$, $\alpha \in (0, \frac{n}{\beta s_G})$, $q \in (0, \infty]$, and $\gamma < \frac{\beta n}{n - \beta \alpha i_G}$ there exists a constant $C = C(\beta, \alpha, n, i_G, s_G, q, |\Omega|) > 0$ such that
        \begin{equation}\label{eq:WaGOm3}
            \left| \left| g(\WaG{f}) \right| \right|_{\Lambda^{\gamma, q}(\Omega)} \leq C||f||_{\Lambda^{\beta, q}(\Omega)}\,
        \end{equation}
        for every $f : \Omega \to \R$ such that $f \in \Lambda^{\beta, q}(\Omega)$.
        \item [(iv)] For any $\alpha \in (0, \frac{n}{ s_G})$, $q \in (0, \infty]$ and $\gamma > \frac{n}{\alpha i_G}$ there exists a constant $C = C(\gamma, \alpha, n, i_G, s_G, q, |\Omega|)$ such that
        \begin{equation}\label{eq:WaGOm4}
            \left| \left| g(\WaG{f}) \right| \right|_{L^{\infty}(\Omega)} \leq C||f||_{\Lambda^{\gamma, q}(\Omega)}\,
        \end{equation}
        for every $f : \Omega \to \R$ such that $f \in \Lambda^{\gamma, q}(\Omega)$.
    \end{enumerate}
\end{prop}

\begin{proof}
We start with the proof of~\eqref{eq:WaGOm1}. Observe that if $t > |\Omega|$, then $\left(g(\WaG{f})\mathds{1}_{\Omega}\right)^*(t) = 0$ and therefore the inequality
\begin{equation*}
    t^{\frac{n - \alpha i_G}{n}}\left(g(\WaG{f})\mathds{1}_{\Omega}\right)^*(t) \leq ||f||_{\Lambda^{1, q}(\Omega)}
\end{equation*}
is satisfied trivially for $t > |\Omega|$. Further, we assume that $t \leq |\Omega|$. Since the function $t \mapsto \frac{g(t)}{t^{i_G - 1}}$ is comparable to a non-decreasing function, we have
\begin{equation*}
    \frac{g \left( |\Omega|^{-\frac{\alpha}{n}}\left(\WaG{f} \right)^*(t) \right)}{|\Omega|^{-\frac{\alpha(i_G - 1)}{n}}} \leq C
    \frac{g \left( t^{-\frac{\alpha}{n}}\left(\WaG{f} \right)^*(t) \right)}{t^{-\frac{\alpha(i_G - 1)}{n}}}\,.
\end{equation*}
Therefore, for some constant $C > 0$ it holds that
\begin{equation*}
    t^{\frac{n - \alpha i_G}{n}}g \left((\WaG{f})^*(t) \right) \leq Ct^{1 - \frac{\alpha}{n}}g(t^{-\frac{\alpha}{n}} (\WaG{f})^{*}(t)).
\end{equation*}
By~\eqref{eq:WaG1} we have
\begin{equation*}
    t^{\frac{n - \alpha i_G}{n}} \left(g\left( \WaG{f}\mathds{1}_{\Omega} \right)\right)^{*}(t) \leq t^{\frac{n - \alpha i_G}{n}}g \left( (\WaG{f})^*(t) \right) \leq C||f||_{\Lambda^{1, q}(\Omega)},
\end{equation*}
which gives us~\eqref{eq:WaGOm1}.

Assertion~\eqref{eq:WaGOm2} can be proved analogically, cf. the proof of Proposition~\ref{prop:WaGest}. By applying it to~\eqref{eq:loremb} we get \eqref{eq:WaGOm3}.

To prove~\eqref{eq:WaGOm4} we observe that since $t \mapsto \frac{g^{-1}(t)}{t^{\frac{1}{i_G-1}}}$ is comparable to a non-increasing function, we have
\begin{equation*}
       t^{\frac{n - \alpha \gamma}{\gamma n (i_G - 1)}}g^{-1}\left( t^{\frac{\alpha}{n}}f^{**}(t)  \right) \leq Cg^{-1}\left( t^{\frac{1}{\gamma}}f^{**}(t) \right)
\end{equation*}
for every $t \leq |\Omega|$. Therefore, by Proposition~\ref{prop:WgOmegaest} it holds that
\begin{align*}
    \left\|\WaG f \right\|_{L^{\infty}(\Omega)} &\leq C\int_{0}^{|\Omega|} t^{\frac{\alpha}{n} - 1}g^{-1}\left( t^{\frac{\alpha}{n}}f^{**}(t) \right)\,dt \leq C\int_{0}^{|\Omega|} t^{\frac{\alpha}{n} - \frac{n - \alpha \gamma}{\gamma n (i_G-1)} - 1}g^{-1} \left( t^{\frac{1}{\gamma}}f^{**}(t) \right)\,dt\\
    &\leq C\int_{0}^{|\Omega|} t^{\frac{\alpha}{n} - \frac{n - \alpha \gamma}{\gamma n (i_G-1)} - 1} g^{-1} \left( \|f\|_{\Lambda^{\gamma, \infty}} \right)\,dt \leq Cg^{-1} \left( \|f\|_{\Lambda^{\gamma, \infty}(\Omega)} \right)\,.
\end{align*}
This gives us that $\left\|g \left(\WaG f \right)\right\|_{L^{\infty}(\Omega)} \leq C\|f\|_{\Lambda^{\gamma, \infty}(\Omega)}$, which is~\eqref{eq:WaGOm4} with $q = \infty$. Using~\eqref{eq:lorinq} gives us~\eqref{eq:WaGOm4} with arbitrary~$q$.
\end{proof}

Using  Hardy inequalities allows formulating more results similar to~\cite[Theorem 3.1]{Ci-pot}.

\begin{prop}\label{prop:WRinfEst}Let  $G\in C(\rp)\cap C^1((0,\infty))$ be an $N$-function such that $G\in\Delta_2\cap\nabla_2$, $\Omega \subset \rn$ be such that $|\Omega| < \infty$, and $\frac{s_G\alpha}{n} \in (0,\min\{1,s_G-1\})$.
\begin{enumerate}
\item [(i)] If $\beta > \max{(1, s_G - 1)}$, and $q \in [1, \infty)$, then there exists a constant $C = C(\alpha, n, i_G, s_G, |\Omega|, \beta, q) > 0$ such that
\begin{equation}\label{eq:Ci1}
    ||\WaG{f} ||_{\Lambda^{\frac{\beta n (s_G - 1)}{n(s_G - 1) - \beta \alpha s_G} , q}(\Omega)} \leq C||g^{-1}(|f|)||_{\Lambda^{\beta,q} (\Omega)}
\end{equation}
for every $f : \Omega \to \R$ such that $g^{-1}(|f|) \in \Lambda^{\beta, q}(\Omega)$.

\item [(ii)] If $s_G \leq 2$ and $q \in (0, 1]$, then there exists a constant $C = C(\alpha, n, i_G, s_G, |\Omega|, q) > 0$ such that
\begin{equation}\label{eq:Ci2}
    ||\WaG{f} ||_{\Lambda^{\frac{n(s_G-1)}{n(s_G - 1) - \alpha s_G} , \infty}(\Omega)} \leq C||g^{-1}(|f|)||_{\Lambda^{1, q}(\Omega)}
\end{equation}
for every $f : \Omega \to \R$ such that $g^{-1}(|f|) \in \Lambda^{1, q}(\Omega)$.

\item [(iii)] There exists a constant $C = C(\alpha, n, i_G, s_G, |\Omega|)$ such that
\begin{equation}\label{eq:Ci3}
    ||\WaG{f} ||_{L^{\infty}(\Omega)} \leq C||g^{-1}(|f|)||_{\Lambda^{\frac{n(s_G-1)}{\alpha s_G} , 1}(\Omega)}
\end{equation}
for every $f : \Omega \to \R$ such that $g^{-1}(|f|)\in \Lambda^{\frac{n(s_G-1)}{\alpha s_G} , 1}(\Omega)$.

\end{enumerate}

\end{prop}

\begin{proof} We start with the proof of auxiliary inequalities. For this we notice that the function $t \mapsto \frac{g^{-1}(t)}{t^{\frac{1}{s_G - 1}}}$ is comparable to a non-decreasing function which enables us to get for $s \leq |\Omega|$ that
\begin{equation*}
   g^{-1}(s^{\frac{\alpha}{n}}f^{**}(s))=g^{-1}\left(|\Omega|^{\frac{\alpha}{n}}\left(\frac{s}{|\Omega|}\right)^{\frac{\alpha}{n}}f^{**}(s)\right) 
   \leq C \left(\frac{s}{|\Omega|}\right)^{\frac{\alpha}{n(s_G - 1)}}{g^{-1}\left(|\Omega|^{\frac{\alpha}{n}}f^{**}(s)\right)} 
   \leq Cs^{\frac{\alpha}{n(s_G - 1)}}g^{-1}(f^{**}(s))\,.
\end{equation*}
 Combining this inequality with Proposition~\ref{prop:WgOmegaest}  gives us that for $t \leq |\Omega|$ we have
\begin{equation}\label{eq:czerwiec16}
    (\WaG f)^*(t) \leq C\int_{\frac{1}{2}t}^{|\Omega|} s^{\frac{\alpha s_G}{n(s_G-1)} - 1}g^{-1}(f^{**}(s))ds\,.
\end{equation}
Let us consider $s_G \leq 2$. We have that the function $g^{-1}$ is comparable to a convex function. Therefore, using Lemma~\ref{lem:psi-rearr} we get
\begin{align*}
    g^{-1}(f^{**}(s)) = g^{-1}\left( \frac{1}{s} \int_{0}^{s} f^{*}(r)dr \right) \leq \frac{C}{s} \int_{0}^{s} g^{-1}(f^*(r))dr = C\left(g^{-1}(|f|) \right)^{**}(s)\,.
\end{align*}
We substitute it into~\eqref{eq:czerwiec16} to obtain that for $s_G \leq 2$ that there exists a constant $C = C(\alpha, n, i_G, s_G, |\Omega|) > 0$ such that for all $t \leq |\Omega|$ and for every measurable $f : \Omega \to \mathbb{R}$ it holds
\begin{align}\label{eq:WgOmegaest}
   (\WaG f)^*(t) &\leq C\int_{\frac{1}{2}t}^{|\Omega|} s^{\frac{\alpha s_G}{n(s_G-1)}}\left(g^{-1}(|f|) \right)^{**}(s)\,\frac{ds}{s}\,.\end{align}
We assume now that $s_G \geq 2$.  We use Lemma~\ref{lem:gconv} and Lemma~\ref{lem:jensen} in the following way
\begin{equation*}
    g^{-1}(f^{**}(s)) \leq Cs^{-\frac{1}{s_G-1}}\int_{0}^{s} r^{\frac{1}{s_G-1}-1}g^{-1}(f^{*}(r))dr\,.
\end{equation*}
By~\eqref{eq:czerwiec16}, for $s_G \geq 2$ there exists a constant $C = C(\alpha, n, i_G, s_G, |\Omega|) > 0$ such that for all $t \leq |\Omega|$ and for every measurable $f : \Omega \to \mathbb{R}$ it holds
\begin{align}
(\WaG f)^*(t) &\leq C\int_{\frac{1}{2}t}^{|\Omega|} s^{\frac{\alpha s_G}{n(s_G - 1)} - \frac{1}{s_G-1}}\int_{0}^{s} r^{\frac{1}{s_G-1}}g^{-1}(f^*(r))\,\frac{dr}{r}\,\frac{ds}{s}\,.\label{eq:WgOmegaest2}
\end{align}  

Now we shall present the proof of~\eqref{eq:Ci1}. Let us denote $\gamma = \frac{\beta n (s_G-1)}{n(s_G-1) - \beta \alpha s_G}$, which is greater than zero by assumptions on $\alpha$. We shall begin with case $s_G < 2$. By~\eqref{eq:WgOmegaest} and Hardy inequality~\eqref{eq:Hardy2} it holds that
\begin{align*}
    ||\WaG f||_{\Lambda^{\gamma, q}(\Omega)}^q 
    &\leq C\int_{0}^{|\Omega|} t^{\frac{q}{\gamma} - 1}\left( \int_{\frac{1}{2}t}^{|\Omega|} s^{\frac{\alpha s_G}{n(s_G-1)} - 1}\left(g^{-1}(|f|) \right)^{**}(s)\,ds \right)^q\,dt\\
    &\leq C \int_{0}^{|\Omega|} t^{\frac{q}{\gamma} + \frac{q\alpha s_G}{n(s_G-1)} - 1} \left( \left(g^{-1}(|f|)\right)^{**}(t) \right)^q\,dt \leq C||g^{-1}(|f|)||_{\Lambda^{\beta, q}(\Omega)}^q,
\end{align*}
where the last inequality holds due to $\frac{1}{\gamma} + \frac{\alpha s_G}{n(s_G-1)} = \frac{1}{\beta} < 1$ and~\eqref{eq:loreq}.
To consider $s_G \geq 2$, we start from~\eqref{eq:WgOmegaest2} and make use of~\eqref{eq:Hardy2} and~\eqref{eq:Hardy1} in the following way.
\begin{align*}
    ||\WaG f||_{\Lambda^{\gamma, q}(\Omega)}^q 
    &\leq C\int_{0}^{|\Omega|} t^{\frac{q}{\gamma} - 1}\left( \int_{\frac{1}{2}t}^{|\Omega|} s^{\frac{\alpha s_G}{n(s_G - 1)} - \frac{1}{s_G-1}} \frac{1}{s}\int_{0}^{s} r^{\frac{1}{s_G-1} - 1}g^{-1}(f^*(r))\,dr\,ds\,  \right)^q\,dt\\
    &\leq C \int_{0}^{|\Omega|} t^{\frac{q}{\gamma} + \frac{q\alpha s_G}{n(s_G-1)} - \frac{q}{s_G - 1} - 1} \left( \int_{0}^{t} r^{\frac{1}{s_G-1} - 1}g^{-1}(f^*(r))\,dr \right)^q\,dt\\ 
    &\leq C \int_{0}^{|\Omega|} t^{\frac{q}{\gamma} + \frac{q\alpha s_G}{n(s_G - 1)} - 1} \left( \left(g^{-1}(|f|)\right)^{*}(t) \right)^q\,dt = C||g^{-1}(|f|)||_{\Lambda^{\beta, q}(\Omega)}^q,
\end{align*}
where the second use of Hardy inequality is permitted due to $\frac{1}{\gamma} + \frac{\alpha s_G}{n(s_G-1)} - \frac{1}{s_G - 1} = \frac{1}{\beta} - \frac{1}{s_G-1} < 0$. \newline

Now we shall look at~\eqref{eq:Ci2}. We denote $\gamma = \frac{n (s_G-1)}{n(s_G-1) - \alpha s_G}$, which is greater than zero due to the assumption on $\alpha$. Using~\eqref{eq:WgOmegaest} leads to 
\begin{align*}
    t^{\frac{1}{\gamma}}\left( \WaG f \right)^*(t) &\leq Ct^{\frac{1}{\gamma}} \int_{\frac{1}{2}t}^{|\Omega|} s^{\frac{\alpha s_G}{n(s_G - 1)} - 2}\int_{0}^{s} g^{-1}(f^*(r))\,dr\,ds\,
    \leq C t^{\frac{1}{\gamma} + \frac{\alpha s_G}{n(s_G-1)} - 1}\|g^{-1}(|f|)\|_{L^1(\Omega)} = C\|g^{-1}(|f|)\|_{L^1(\Omega)}\,.
\end{align*}
The inequality~\eqref{eq:Ci2} is therefore proven for $q = 1$ and using~\eqref{eq:lorinq} leads to~\eqref{eq:Ci2} itself.  \newline 

Now we shall prove~\eqref{eq:Ci3}. In case of $s_G < 2$, we observe that since $\alpha < \frac{n(s_G-1)}{s_G}$ and we have~\eqref{eq:WgOmegaest}, it holds that
\begin{equation*}
    || \WaG f ||_{L^{\infty}(\Omega)} \leq C\int_{0}^{|\Omega|} t^{\frac{\alpha s_G}{n(s_G-1)} - 1}(g^{-1}(|f|))^{**}(t)\,dt \leq C||g^{-1}(|f|)||_{\Lambda^{\frac{n(s_G-1)}{\alpha s_G}, 1}(\Omega)}.
\end{equation*}
To achieve a similar result for $s_G \geq 2$, we use~\eqref{eq:WgOmegaest2} in conjunction with Hardy inequality~\eqref{eq:Hardy1} and get
\begin{align*}
    || \WaG f ||_{L^{\infty}(\Omega)} 
    &\leq C\int_{0}^{|\Omega|} t^{\frac{\alpha s_G}{n(s_G - 1)} - \frac{1}{s_G-1}} \frac{1}{t}\int_{0}^{t} r^{\frac{1}{s_G-1} - 1}g^{-1}(f^*(r))\,dr\,ds\,\\
    &\leq C \int_{0}^{|\Omega|} t^{\frac{\alpha s_G}{n(s_G-1)} - 1}g^{-1}(f^{*}(t))\,ds = C||g^{-1}(|f|)||_{\Lambda^{\frac{n(s_G-1)}{\alpha s_G}, 1}(\Omega)}.
\end{align*}
Hence, inequality~\eqref{eq:Ci3} is proven.
\end{proof}

 It is natural to expect a result on boundedness of $\WaG$ in the Orlicz scale.  For an $N$-function $A$, by an Orlicz space ${L}^A(\Omega)$  we understand the space of measurable functions $f:\Omega\to\R$ such that there exists $\lambda>0$ making $\int_{\Omega}A\left( \tfrac{1}{\lambda} |f|\right)\,dx<\infty$, endowed with the Luxemburg norm 
\[||f||_{L^A(\Omega)}:=\inf\left\{\lambda>0:\ \ \int_{\Omega}A\left( \tfrac{1}{\lambda} |f|\right)\,dx\leq 1\right\}\,.\]
Let us present an application of~\cite[Theorem 3.4]{Ci-pot}, which is proven via interpolation arguments,  in order to infer that $g\left(\WaG\right):L^A\to L^B$ for compatible $A$ and $B$.


\begin{coro}\label{coro:Orlicz}
Let $\alpha < \frac{n}{s_G}$. Suppose $A$ and $B$ are $N$-functions and let $\Omega \subset \rn$ be such that $|\Omega| < \infty$. For $\sigma>\frac{n}{\alpha i_G}$ we define
\begin{equation}\label{eq:EFdef}
    E(t) := \left( \int_{0}^{t} \left( \frac{r}{A(r)}  \right)^{\frac{1}{\sigma-1}} \,dr \right)^{\frac{\sigma-1}{\sigma}}\qquad\text{and}\qquad
    F(t) := \left( \int_{0}^{t} \frac{B(r)}{r^{1 + \frac{\sigma}{\sigma-1}}}\,dr \right)^{\frac{\sigma-1}{\sigma}}\,,
\end{equation}
and suppose that $E$ is finite-valued for $t > 0$. Assume further that there exists a constant $\delta > 0$ such that for $t > 0$ it holds that
\begin{equation}\label{compatibility}
    F\left( \frac{E(t)}{\delta} \right) \leq \delta \frac{A(t)}{t}\,.
\end{equation}
Then there exists a constant $C = C(\alpha, n, i_G, s_G, |\Omega|, \delta)$ such that
\begin{equation*}
    ||g(\WaG f)||_{L^B(\Omega)} \leq C||f||_{L^A(\Omega)}
\end{equation*}
for every $f : \Omega \to \R$ satisfying $f \in L^A(\Omega)$.
\end{coro}
\begin{proof} Despite the operator $\mathcal{T}(f):=g(\WaG f)$ is not a quasilinear operator as it does not satisfy $|\lambda|\mathcal{T}(f)=\mathcal{T}(\lambda f)$, we note that the proof of~\cite[Theorem 3.4]{Ci-pot} does not use this property. In order to apply this result, we need to show that $\mathcal{T}:L^1\to \Lambda^{\frac{\sigma}{\sigma-1},\infty}$ and $\mathcal{T}: \Lambda^{\sigma,1}\to L^\infty$. Boundedness of $\mathcal{T}$ between $\Lambda^{\sigma, 1}$ and $L^{\infty}$ comes from~\eqref{eq:WaGOm4}, while boundedness between $L^1$ and $\Lambda^{\frac{\sigma}{\sigma - 1}, \infty}$ comes from~\eqref{eq:WaGOm1} and~\eqref{eq:loremb}. We are in position to apply~\cite[Theorem 3.4]{Ci-pot} to conclude.
\end{proof}

\subsection{Regularity of solutions to elliptic equations in reflexive Orlicz spaces}\label{sec:regularity}

Potential estimates are known to be precise tools in the regularity theory, see \cite{KuMiguide} and references therein. In this section we follow the ideas of \cite{Ci-pot} and \cite{CiSch} dealing with the $p$-Laplace problems in order to transfer regularity from the datum to the solution also in the case of problems involving Orlicz growth being of broad interest recently, see e.g.,~\cite{Baroni-Riesz,MBIC,CiMa,IC-gradest,IC-pocket} and references therein. We shall employ known estimates of \cite{Baroni-Riesz,CGZG-Wolff,CYZG-Wolff,ACCZG}.\newline

Within this section, we concentrate on the regularity of local weak solutions (when data is regular enough) and solutions obtained in an approximation procedure (when data is too poorly integrable for weak solutions to exist). The equation we deal with has a form~\eqref{eq-Or}, which is studied under the following regime.\newline

\underline{Assumption {\bf (A)}}. Given a bounded open set $\Omega\subset\rn$, $n\geq 2$, let us consider a~Carath\'eodory's function $\opA:\Omega\times \rn\to\rn$ (measurable with respect to the first variable and continuous with respect to the second one), which is strictly monotone, that is 
\[(\opA(x,\xi)-\opA(x,\eta))\cdot(\xi-\eta)>0\qquad\text{for every } \ \xi\neq\eta.\]
We suppose that $\opA$ satisfies growth and coercivity conditions expressed by the means of a doubling $N$-function $G\in C^1 {((0,\infty))} \cap C([0, \infty))$ i.e., a nonnegative, increasing, and convex function such that $G\in\Delta_2\cap\nabla_2$. Namely, we assume that
\begin{equation}
\label{ass-op}\begin{cases}
c_1^\opA G(|\xi|)\leq \opA(x,\xi)\cdot\xi,\\
|\opA(x,\xi)|\leq c_2^\opA g(|\xi|),
\end{cases}
\end{equation}
where $g$ is the derivative of $G$ and $c_1^\opA,c_2^\opA>0$ are absolute constants. We collect all parameters of the problem as $ \data=\data(i_G,s_G,c_1^\opA,c_2^\opA),$ where $i_G$ and $s_G$ describe the growth of $G$, cf.~\eqref{iG-sG}.\newline
 
 We define the Orlicz--Sobolev space  $W^{1,G}(\Omega)$  as follows
\begin{equation*} 
W^{1,G}(\Omega)=\big\{f\in W^{1,1}_{loc}(\Omega):\ \ |f|,|D f|\in L^G(\Omega)\big\},
\end{equation*}where the gradient is understood in the distributional sense, endowed with the norm
\[
\|f\|_{W^{1,G}(\Omega)}=\inf\bigg\{\lambda>0 :\ \    \int_{\Omega}G\left( \tfrac{1}{\lambda} |f|\right)+ G\left( \tfrac{1}{\lambda} |Df|\right)\,dx\leq 1\bigg\} 
\]
and  by $W_0^{1,G}(\Omega)$ we denote a closure of $C_0^\infty(\Omega)$ under the above norm. 
Properties of these function spaces depend only on the behaviour of $G$ only up to equivalence and far from the origin.\newline

A function $u\in W^{1,G}_{loc} (\Omega)$ is called a local weak solution to equation \eqref{eq-Or} under Assumption {\bf (A)} if
\[\int_{\Omega'}\opA(x, Du)\cdot D\vp\,dx=\int_{\Omega'} f\vp\,dx\quad\text{
for every open set $\Omega'\Subset\Omega$ and every $\vp\in W^{1,G}_0(\Omega')$.}\]

Recall that we define the Wolff-type potential $\cW_{1, G}$ in~\eqref{wolff-potential}. Following~\eqref{eq:WR-psi} and choosing $\psi(t)=g^{-1}(t)$ we consider the truncated potential given by
\begin{equation}\label{eq:WR} \cW_{1, G}^{\blue{R}} f(x) = \int_{0}^{R} g^{-1}\left(r^{1-n}{\int_{B(x,r)}|f(y)|\,dy}\right)\,dr. \end{equation}

By \cite[Theorem~2]{CGZG-Wolff} we know that an $\opA$-superharmonic function $u$ generated by a nonnegative $f\in L^1_{loc}$  satisfies~\eqref{est-wolff-G}. In this situation $\opA$-superharmonic functions are defined via the comparison principle with local continuous solutions to the homogeneous equation $-\dv \opA(x,Du)=0$. If $f\geq 0$ is regular enough (i.e., $f\in (W^{1,G}_0(\Omega))^*$), the $\opA$-superharmonic function $u$ generated by $f$ is a local weak solution. In turn, local weak solutions are controlled on small balls by $\cW_{1,G}^R$ via~\eqref{est-wolff-G}. Let us make use of Theorem~\ref{theo:equivalence-G} to find more consequences of the pointwise estimate~\eqref{est-wolff-G}  for the solution $u$.

\begin{coro}\label{coro:hu-Y-X} Suppose $G,n,\opA, \Omega$ satisfy Assumption {\bf (A)} and $s_G < n$. Let $X(\Omega)$ and $Y(\Omega)$ be quasi-normed rearrangement invariant spaces and let $h_X, h_Y : \rp \to \rp$ be non-decreasing, left-continuous functions such that $h_X, h_Y \in \Delta_2$. Assume further that for every \blue{$a \in (0, |\Omega|)$} there exists a constant $c>0$ such that for every nonnegative function $\phi \in \overline{X}(0,\blue{a})$ it holds
    \begin{equation}\label{eq:ass-XY-est} \left\| h_Y \left(\int_{\frac{1}{2}t}^{\blue{a}} s^{\frac{1}{n}-1}g^{-1}\left(s^{\frac{1}{n}-1}\int_{0}^{s}\phi(y)\,dy\right)ds \right) \right\|_{\overline{Y}(0,\blue{a})} \leq c ||h_X(\phi)||_{\overline{X}(0,\blue{a})}. \end{equation}
Then a local weak solution $u\in W^{1,G}_{loc}(\Omega)$ to~\eqref{eq-Or} with nonnegative $f\in \blue{X_{loc}(\Omega)}$ satisfies $u\in{Y_{loc} (\Omega)}$.\\ Moreover, for every $\Omega''\Subset\Omega'\Subset\Omega$ it holds that
\begin{equation}
    \label{hu-Y-X}
    \|h_Y(|u|)\|_{Y(\Omega'')}\leq C\left(h_Y\left(g^{-1}\left(\|g(|u|)\|_{L^1(\Omega')}\right)\right)+{\rm diam}\, \Omega+\|h_X(|f|)\|_{X(\Omega')}\right). \end{equation}\end{coro}
\begin{proof}

Let us take any $R > 0$ such that~\eqref{est-wolff-G} holds true and satisfying $R < \tfrac12 \text{dist}(\partial \Omega', \partial \Omega'')$, so that $B(x, R) \subseteq \Omega'$ for all $x \in \Omega''$. Observe that for such $x$ we have $\cW_{1,G}^R f(x) = \cW_{1,G}^R (f\mathds{1}_{\Omega'})(x)$. Moreover, we can estimate
\begin{equation*}
\inf_{B(x, R)} |u(y)| = g^{-1}\left(\inf_{B(x, R)} g(|u(y)|) \right) \leq g^{-1}\left(\frac{1}{\omega_n R^n} \int_{B(x, R)} g(|u(y)|)\,dy \right) \leq g^{-1} \left(\frac{1}{\omega_n R^n} \|g(|u|)\|_{L^1(\Omega')} \right).
\end{equation*}
Therefore, by~\eqref{est-wolff-G} we have an estimate of the form
\begin{equation*}
    |u(x)| \leq C\left( g^{-1}\left(\|g(|u|)\|_{L^1(\Omega')}\right) + {\rm diam}\, \Omega + \cW_{1,G}^R (f\mathds{1}_{\Omega'})(x) \right)\,.
\end{equation*}
Applying $h_Y$ and using its doubling properties give us that
\begin{equation*}
    h_Y(|u(x)|) \leq C\left( h_Y \left(g^{-1}\left(\|g(|u|)\|_{L^1(\Omega')}\right) \right) + {\rm diam}\, \Omega + h_Y \left(\cW_{1,G}^R (f\mathds{1}_{\Omega'})(x) \right) \right)\,.
\end{equation*}
%
%
%
By using Theorem~\ref{theo:equivalence-G} and~\eqref{eq:ass-XY-est} we can estimate
\begin{align*}
||h_Y \left(\cW_{1,G}^R (f\mathds{1}_{\Omega'}) \right)||_{Y(\Omega'')} &\leq ||h_Y \left(\cW_{1,G}^R (f\mathds{1}_{\Omega'}) \right)||_{Y(\Omega')} \leq C||h_X(|f|)||_{X(\Omega')}\,.
\end{align*}
Therefore, 
we have
\begin{equation}\label{eq:lipiec13}
    ||h_Y(|u|)||_{Y(\Omega'')} \leq C\left( h_Y \left(g^{-1}\left(\|g(|u|)\|_{L^1(\Omega')}\right) \right) + {\rm diam}\, \Omega + \|h_X(|f|)\|_{X(\Omega')} \right)\,.
\end{equation}
Note that by~\cite[Lemma 4.1]{CiMa} it holds that $g(|u|) \in L^1_{loc}(\Omega)$, which means that the right-hand side of~\eqref{eq:lipiec13} is finite.
\end{proof}

We can use Corollary~\ref{coro:hu-Y-X} and results from the previous section to obtain regularity results for $u$. In fact, we apply  Proposition~\ref{prop:WRinfEst}, as well as Proposition~\ref{prop:gWest}  to get what follows. 
\begin{theo}\label{theo:solution1}
Suppose $G,n,\opA, \Omega$ satisfy Assumption {\bf (A)}, $s_G < n$ and $u$ is a local weak solution to~\eqref{eq-Or} with nonnegative data. Then\begin{itemize}
    \item [(i)] if $q \in (0, 1]$ and $f \in \Lambda^{1, q}_{loc}(\Omega)$, then $g(|u|) \in \Lambda^{\frac{n}{n - i_G}, \infty}_{loc}(\Omega)$;
    \item [(ii)] if $\beta > 1$, $q \in [1, \infty]$ and $f\in \Lambda^{\beta, q}_{loc}(\Omega)$, then $g(|u|)\in \Lambda^{\frac{\beta n}{n - \beta i_G} , \infty}_{loc}(\Omega)$; 
    \item [(iii)] if $q \in (1, \infty]$, $\gamma > \frac{n}{i_G}$ and $f \in \Lambda^{\gamma, q}_{loc}(\Omega)$, then $u \in L^{\infty}_{loc}(\Omega)$;
    \item [(iv)] if $\frac{n}{n-1} < s_G \leq 2$, $q \in (0, 1]$ and $g^{-1}(|f|) \in \Lambda^{1, q}_{loc}(\Omega)$, then $u \in \Lambda^{\frac{n(s_G-1)}{n(s_G-1) - s_G}, \infty}_{loc}(\Omega)$;
    \item [(v)] if $\beta \in \left(\max(1, s_G - 1), \frac{n(s_G-1)}{s_G} \right)$, $q \in [1, \infty)$ and $g^{-1}(|f|) \in \Lambda^{\beta, q}_{loc}(\Omega)$, then $u \in \Lambda^{\frac{\beta n(s_G-1)}{n(s_G-1)-\beta s_G} , q}_{loc}(\Omega)$;
    \item [(vi)] if $\frac{n}{n-1} < s_G$ and $g^{-1}(|f|) \in \Lambda^{\frac{n(s_G-1)}{s_G}, 1}_{loc}(\Omega)$, then $u \in L^{\infty}_{loc}(\Omega)$\,.
\end{itemize}
\end{theo}

In particular, in the view of~\eqref{eq:loremb} and Lemma~\ref{lem:zyginv}, for $G$ being a Zygmund function we have what follows. 

\begin{coro}\label{coro:solution2}
Let $1 < p < n$, $\gamma \in \R$. Suppose $n,\opA, \Omega$ satisfy Assumption {\bf (A)} with $G(t) = t^p\log^{\gamma}(s + t)$ for sufficiently large $s$, and $u$ is a local weak solution to~\eqref{eq-Or}. Then\begin{itemize}
    \item [(i)] if $q \in (0, 1]$ and $f \in \Lambda^{1, q}_{loc}(\Omega)$, then $|u|^{p-1}\log^{\gamma}(s+|u|) \in \Lambda^{\delta, \infty}_{loc}(\Omega)$ for every $0< \delta < \frac{n}{n-p}$;
    \item [(ii)] if $\beta > 1$, $q \in [1, \infty]$ and $f\in \Lambda^{\beta, q}_{loc}(\Omega)$, then $|u|^{p-1}\log^{\gamma}(s+|u|)\in \Lambda^{\delta , \infty}_{loc}(\Omega)$ for every $0 < \delta < \frac{\beta n}{n - \beta p}$;
    \item [(iii)] if $q \in (1, \infty]$, $\gamma > \frac{n}{p}$ and $f \in \Lambda^{\gamma, q}_{loc}(\Omega)$, then $u \in L^{\infty}_{loc}(\Omega)$;
    \item [(iv)] if $\frac{n}{n-1} < p < 2$, $q \in (0, 1]$ and $|f|^{\frac{1}{p-1}}
    \left(\log{(s + |f|)} \right)^{-\frac{\gamma}{p-1}} \in \Lambda^{1, q}_{loc}(\Omega)$, then $u \in \Lambda^{\delta, \infty}_{loc}(\Omega)$ for every $0 < \delta < \frac{n(p-1)}{n(p-1)-p}$;
    \item [(v)] if $\beta \in \left(\max\{1, p -1\}, \frac{n}{p}(p-1) \right)$, $q \in [1, \infty)$ and $|f|^{\frac{1}{p-1}}
    \left(\log{(s + |f|)} \right)^{-\frac{\gamma}{p-1}} \in \Lambda^{\beta, q}_{loc}(\Omega)$, then \mbox{$u \in \Lambda^{\delta , q}_{loc}(\Omega)$} for every $0 < \delta < \frac{\beta n(p-1)}{n(p-1)-\beta p}$;
    \item [(vi)] if $p > \frac{n}{n-1}$ and $|f|^{\frac{1}{p-1}}
    \left(\log{(s + |f|)} \right)^{-\frac{\gamma}{p-1}} \in \Lambda^{\delta, 1}_{loc}(\Omega)$ for some $\delta > \frac{n}{p}(p-1)$, then $u \in L^{\infty}_{loc}(\Omega)$\,.
\end{itemize}
\end{coro}
Results similar to Corollary~\ref{coro:solution2} may be obtained for problems governed by $G(t) = t^p\left(\log \log (s+t)\right)^{\gamma}$ by using Lemma~\ref{lem:zyginv2}. We left the exact formulation to the reader. 

 \blue{For results related to {\it (iii)} for problems with more general structure involving also explicit $u$-dependence of an operator see e.g. \cite{Marcellini1}.}

\medskip

Using Corollary~\ref{coro:Orlicz} we can also infer the following result.

\begin{theo}\label{thm:orlicz-result} Suppose $G,n,\opA, \Omega$ satisfy Assumption {\bf (A)}, and $u$ is a local weak solution to~\eqref{eq-Or} with nonnegative data. 
For $\sigma>\frac{n}{i_G}$ let $E$ and $F$ be defined as in~\eqref{eq:EFdef}. Assume further that $E$ is finite-valued for $t > 0$ and there exists a constant $\delta > 0$ such that for $t > 0$ such that~\eqref{compatibility} holds true. If $|f| \in L^{A}_{loc}(\Omega)$, then $g(|u|) \in L^B_{loc}(\Omega)$.
\end{theo}
We refer to \cite{CiSch} for the potential estimates for solutions to $p$-Laplace systems of equations with regular data in divergence form and more methods of using boundedness of $\Wap$, $p>1$, to infer regularity of solutions. The results are formulated in the general scales of Morrey-type and Campanato-type. This approach seems to  be possible to adapt to the vectorial and the scalar case with Orlicz growth and measurable coefficients with the use of the potential $\WaG$. However, it is not direct and cannot be used to stress applicability of Theorems~\ref{theo:boundedness-psi} and~\ref{theo:equivalence}.

Let us mention other results that can be inferred for vectorial problems. For a result related to the upper bound from~\eqref{est-wolff-G} for solutions to systems with Orlicz and superquadratic growth and under more typical assumptions like quasidiagonal structure of the system, we refer to \cite{CYZG-Wolff} proven with the method of \cite{KuMi2018jems}. 
\begin{rem}\rm \label{rem:vectorial} For a local weak solution $\mathbf{u}:\Omega\to\r^m$ to system $-\mathbf{\rm div} \opA(x,D\mathbf{u})=\mathbf{f}$ with a datum $\mathbf{f}\in (L^1(\Omega))^m$, \cite[Theorem~2.1]{CYZG-Wolff} yields that for all sufficiently small $R>0$ and $B(x,R)\Subset B(x,2R)\Subset\Omega$ it holds that
\[ |\mathbf{u}(x )|\leq C_\cW\left(\cW_{1,G}^R(|\mathbf{f}|) (x)+\frac{1}{|B(x,R)|}\int_{B(x,R)} |\mathbf{u}(y)|dy\right).\]
Consequently, $|\mathbf{u}|$ shares the regularity consequences from Theorem~\ref{theo:solution1} and can be illustrated by a related counterpart of Corollary~\ref{coro:solution2}. In particular, this result concerns solutions to vectorial $p$-Laplace problems studied in~\cite{KuMi2018jems}.
\end{rem}

\begin{rem}
\label{rem:SOLA}\rm
In fact, Corollary~\ref{coro:hu-Y-X} requires only $g(|u|)\in L^{1}_{loc}(\Omega)$. In turn, one can infer from estimate~\eqref{hu-Y-X} all consequences of Corollary~\ref{coro:hu-Y-X}, Theorem~\ref{theo:solution1}, Remark~\ref{rem:vectorial} for a more general notion of solutions e.g., for SOLA (solutions obtained as a limit of approximation) introduced in~\cite{bgSOLA}, approximable solutions~\cite{CiMa,ACCZG}, or $\opA$-superharmonic functions~\cite{HeKiMa,CGZG-Wolff}. This passage to the limit in the Orlicz setting in the case of SOLA is described in~\cite{Baroni-Riesz}, while for $\opA$-superharmonic functions in~\cite{CGZG-Wolff}. Very weak solutions to a system from Remark~\ref{rem:vectorial} are proven to exist in~\cite{CYZG2}. 
\end{rem}

Let us note that gradients of solutions to~\eqref{eq-Or} are also known to be controlled by a potential, see~\cite{Ci-pot,KuMiguide}. In fact, potential estimates on $|Du|$, where $u$ is a SOLA to \eqref{plap}, were provided in the terms of the Wolff potential in \cite{DuMiJFA2010
} and then improved to estimates of $|Du|^{p-1}$ by a linear Riesz potential $\cI_1$ in~\cite{KuMiARMA2013}. Note that due to~\cite[(26)]{KuMiguide} for $p\geq 2$ it holds that $\cI_1^R f\lesssim (\cW_{\frac{1}{p},p}^{2R}f)^{p-1}$. For a natural subclass of \eqref{eq-Or} in \cite{Baroni-Riesz} it is proven that\begin{equation}\label{est-baroni}
    g(|Du(x)|) \leq C\left(\cI_1^{2R} f(x) +  g\left(\frac{1}{|B(x,R)|}\int_{B(x,R)}|Du|\, dx\right)\right),\qquad\text{where }\ C=C(\data),
\end{equation}
and $\cI_1^R f$ is the truncated Riesz potential, that is, $\cI_1^R f(x):=\int_0^R\left(r^{-n}\int_{B(x,r)}|f|\,dy\right)\,dr\,.$ Note that for SOLA it is known that $g(|Du|)\in L_{loc}^1(\Omega)$.  Boundedness of Riesz potential is studied between Lorentz spaces, $L^1$ and Marcinkiewicz spaces, as well as Lorentz space and $L^\infty$ are established in \cite[Proposition 2.8]{Ci-pot}, see also \cite{Cianchi-interp,Cianchi-strong}. Observe that truncated Riesz potential can be estimated using maximal function, i.e.,
\begin{equation*}
    \cI_1^R f(x) \leq \int_{0}^{R} M_1f(x)\,dr = R M_1f(x).
\end{equation*}
Therefore, $\cI_1^R$ shares the same estimates as the maximal operator. In turn, using Lemma~\ref{lem:maximal-function} and~\cite[Proposition 2.8]{Ci-pot} as well, we can formulate the counterpart of Theorem~\ref{theo:solution1} for gradients of solutions.

\begin{theo}\label{theo:gradient}
Suppose $G,n,\opA$ satisfy Assumption {\bf (A)}.  Let us additionally assume that $G$ is superquadratic (i.e., $i_G\geq 2$), $\opA(x,\xi)=\opA(\xi)$,  and $u$ is a local weak solution to~\eqref{eq-Or}.  Then the following implications hold true.\begin{itemize}
    \item[(i)] If $1 < \sigma <n$
and $0 <\beta\leq\infty$ and $f\in\Lambda^{\sigma,\beta}_{loc}(\Omega)$, then $g(|Du|)\in \Lambda^{\frac{\sigma n}{n-\sigma},\beta}_{loc}(\Omega)$;
   \item[(ii)]  If $f\in L^1_{loc} (\Omega)$, then $g(|Du|)\in \Lambda^{\frac{n}{n-1},\infty}_{loc}(\Omega)$;
   \item[(iii)] If locally $f\in L\log L (\Omega)$, then $g(|Du|)\in L_{loc}^{\frac{n}{n-1}}(\Omega)$;
   \item[(iv)] If $1 < q < \theta \leq n$ and locally $f \in L^{q}_{\theta}(\Omega)$, then locally $g(|Du|) \in L_{\theta}^{\frac{\theta q}{\theta - q}}(\Omega)$;
   \item[(v)] If $1 < q < \theta \leq n$, $s \in (0, \infty)$ 
 and locally $f \in \Lambda^{q, s}_{\theta}(\Omega)$, then locally $g(|Du|) \in \Lambda^{\frac{\theta q}{\theta - q}, \frac{\theta s}{\theta - q}}(\Omega)$\,.
 \end{itemize}
\end{theo}
A version of the above theorem extended to systems with Orlicz growth similar to described in Remark~\ref{rem:vectorial}  is provided in~\cite{CKW}. The results {\it (i)--(ii)} above are closely related to~\cite[Theorem~1]{IC-gradest}. Theorem~\ref{theo:gradient} provides the regularity for a larger set of parameters. On the other hand, \cite[Theorem~1]{IC-gradest} allows for analysis of the operator with measurable coefficients, i.e., $\opA:\Omega\times\rn\to\rn$ being a Carath\'eodory's function.\newline

Given the current interest in Lipschitz bounds, see e.g.,~\cite{BeckMingione,DiMM}, let us formulate one more consequence of~\eqref{est-baroni} and boundedness of $\cI_{1} : \Lambda^{n, 1}(\rn) \to L^{\infty}(\rn)$~\cite[Proposition 2.8]{Ci-pot}. Note that the key missing ingredient to get rid of the technical assumptions $i_G\geq 2$ and $\opA(x,\xi)=\opA(\xi)$ is to extend~\eqref{est-baroni} coming from \cite{Baroni-Riesz} to this case.
\begin{coro}[{Lipschitz  bound}]\label{coro:Lip}
Suppose $G,n,\opA$ satisfy Assumption {\bf (A)}.  Let us additionally assume that $G$ is superquadratic (i.e., $i_G\geq 2$), $\opA(x,\xi)=\opA(\xi),$ and $u$ is a local weak solution to~\eqref{eq-Or}. 
 If $f\in \Lambda^{n,1}_{loc}(\Omega)$, then $g(|Du|)\in L^\infty_{loc}(\Omega)$. Moreover, for every $\Omega''\Subset\Omega'\Subset\Omega$ 
 we have\begin{equation*} 
    \|g(|Du|)\|_{L^\infty(\Omega'')} \leq C\left(\|f\|_{\Lambda^{n,1}(\Omega')} +  g\left(\|\,|Du|\,\|_{L^1(\Omega')}\right)\right),\qquad\text{where }\ C=C(\data)\,.
\end{equation*}
\end{coro}
\noindent As in the case of a solution $u$ described in Remark~\ref{rem:SOLA}, estimate~\eqref{est-baroni} for gradients work also for gradient of SOLA, see~\cite{Baroni-Riesz}. In turn, the claims of Theorem~\ref{theo:gradient} and Corollary~\ref{coro:Lip} are true for $u$ being SOLA. Similar bound for solutions to vectorial problems is proven in~\cite{CKW}. Let us stress, however, that \cite[Theorem~1.15]{BeckMingione} provides also this kind of Lipschitz estimate for energy solutions to vectorial problems for more general case, i.e., without assuming that the growth of the operator is superquadratic. \blue{For related results for problems with more general structure involving also explicit $u$-dependence of an operator see e.g. \cite{Marcellini2,Marcellini3}.}

\section{New regularity results for solutions to problems in non-reflexive Orlicz spaces}\label{sec:non-doub}

We point out that the estimates of Theorem~\ref{theo:equivalence} can be used also in order to infer regularity to {nondoubling} problems anisotropic considered e.g. in~\cite{Ci-sym,ACCZG}.

Let us consider the problem $-\dv\opA(x,Du)=f$, where the growth of a monotone Carath\'eodory's function $\opA:\Omega\times\rn\to\rn$ is governed by an anisotropic $N$-function
$\Phi:\rn\to\rp$ via conditions $\opA(x,\xi)\cdot \xi \geq \Phi(\xi)$ and $\widetilde{\Phi}(c_\Phi \opA(x, \xi))\leq  \Phi(\xi) +h(x)$ holding for every  $\xi\in\rn$, a.e. $x\in\Omega$, a positive constant $c_{\Phi}$, and nonnegative $h\in L^1(\Omega)$, where $\widetilde{\Phi}$ denotes the Young conjugate of $\Phi$. For any weak solution $u$, symmetrization methods of \cite{Ci-sym} lead to 
\[u^*(|x|)\leq \frac{1}{n \omega_n^{1/n}}\int _{|B(0,|x|)|}^{|\Omega|} s^{\frac{1}{n}-1}\Psi_\Diamond^{-1}\bigg(\frac{
s^{\frac{1}{n}}}{n \omega_n^{1/n}}f^{**}(s)\bigg)\, ds=:I.\]
In this formula $\Psi_\Diamond$ is equivalent (up to absolute constants) to $\Psi_\circ(t)=\Phi_\circ(t)/t$ where function $ \xi \mapsto \Phi _\circ(|\xi|)$ can be regarded as a kind of  \lq\lq average in measure\rq\rq  \, of $\Phi$, i.e., $\Phi_\circ : [0, \infty ) \to [0, \infty)$ obeys $|\{\xi \in \rn: \Phi_\circ (|\xi|) \leq t\}|  =|\{\xi \in \rn: \Phi
(\xi)\leq t\}|$ for $t\geq 0$. We assume that $\Psi_\Diamond(s)\to\infty$ as $s\to\infty$. It is important to stress that there is no $\Delta_2$-condition imposed on $\Phi$. Using the same arguments as in the proofs of Theorem~\ref{theo:equivalence}   and Theorem~\ref{theo:equivalence-operators}   combined with Remark~\ref{rem:exTx}  one can infer the following result.

\begin{coro}\label{coro:aniso-hu-Y-X} Suppose $\Phi$ and $\opA$ are as above and there exists a weak solution $u$ to $-\dv\opA(x,Du)=f\in L^1_{loc}(\Omega)$. Let $X(\Omega)$ and $Y(\Omega)$ be quasi-normed rearrangement invariant spaces. 

\begin{enumerate}
    \item [(i)] If there exists a constant $c>0$ such that for every nonnegative function $\phi \in \overline{X}(0,\infty)$ it holds
    \begin{equation*} \left\| \int_{|B(0,\cdot)|}^{|\Omega|} s^{\frac{1}{n}-1}\Psi_\Diamond^{-1}\left(\frac{
s^{\frac{1}{n}-1}}{n \omega_n^{1/n}} \int_{0}^{s}\phi(y)\,dy\right)ds \right\|_{\overline{Y}(0,|\Omega|)} \leq c ||\phi||_{\overline{X}(0,|\Omega|)}\,, \end{equation*}
then a weak solution $u$ to~\eqref{eq-Or} with $f\in X(\Omega)$ satisfies $u \in{Y (\Omega)}$. Moreover, \[ \|u\|_{Y(\Omega)}\leq C\|f\|_{X(\Omega)}. \]
\item[(ii)] Assume that $H_X, H_Y : \rp \to \rp$ are non-decreasing functions and $\alpha,\beta,\sigma,\vr\in\R$. If for some $k_1 > 0$ there exists a~constant $c_1>0$ such that for every nonnegative function $\phi \in \overline{X}(0, |\Omega|)$ it holds
    \begin{equation*} \int_0^\infty t^\alpha H_Y\left[k_1 t^\beta \int_{|B(0,t)|}^{|\Omega|} s^{\frac{1}{n}-1}\Psi_\Diamond^{-1}\left(\frac{
s^{\frac{1}{n}-1}}{n \omega_n^{1/n}}\int_{0}^{s}\phi(y)dy\right)ds\right] \,dt\leq \int_0^\infty t^\sigma H_X\left(c_1 t^\vr \phi(t)\right)\,dt\,, \end{equation*} 
    then there exist $k_2>0$ and $c_2 > 0$ such that a weak solution $u$ to~\eqref{eq-Or} with $f \in L^1_{loc}(\Omega)$ satisfies
    \[\int_{0}^{\infty} t^{\alpha}H_{Y} \left(k_2t^{\beta}u^*(t) \right)\,dt \leq  \int_{0}^{\infty} t^{\sigma}H_X \left( c_2t^{\varrho}f^*(t) \right)\,dt\,.\]
 \end{enumerate}
\end{coro}
By arguments from Remark~\ref{rem:SOLA} one can pass to the limit in the case of Corollary~\ref{coro:aniso-hu-Y-X} for approximable solutions from~\cite{ACCZG}.

Note that the results above are specializing to the isotropic situation with operators exposing Orlicz and non-doubling growth, for which this result is already new.

\section{Appendix}

\subsection{Relation between Havin--Maz'ya and Wolff potentials}
\label{ssec:rem:V=W}
Let us consider a general form of $\cV_{\alpha,p}$ from~\eqref{havin-mazya-potential}. Namely, for any non-decreasing function $\psi:\rp\to\rp$ we define
\begin{equation}
    \label{havin-mazya-potential-psi}\bV_{\alpha,\psi}f(x):= \cIa\psi\big(\cIa f \big) (x)\,
\end{equation}
acting on measurable functions $f$ such that~\eqref{skoncz-wart} holds true. Note that for $\psi(s)=s^\frac{1}{p-1}$ we have $\cV_{\alpha,p}=\bV_{\alpha,\psi}$. There exist $c_1 = c_1(\alpha, n) > 0$ and  $c_2 = c_2(\alpha, n) > 0$ such that the pointwise estimate \[c_1{\Waps}(c_2f)(x) \leq \bV_{\alpha, \psi}f(x)\] holds true. Indeed,
\begin{align*}
\bV_{\alpha, \psi}f(x) &= \int_{0}^{\infty} r^{\alpha - n - 1} \int_{B(x, r)}\psi \left( \int_{0}^{\infty} s^{\alpha - n - 1}\int_{B(y, s)}|f(z)|\,dz\,ds \right)\,dy\,dr
\\&\geq \int_{0}^{\infty} r^{\alpha - n - 1} \int_{B(x, r)}\psi \left( \int_{2r}^{\infty} s^{\alpha - n - 1}\int_{B(x, r)}|f(z)|\,dz\,ds \right)\,dy\,dr
\\&= \omega_n\int_{0}^{\infty} r^{\alpha - 1} \psi \left( \tfrac{2^{\alpha - n}}{n - \alpha} \cdot r^{\alpha - n}\int_{B(x, r)}|f(z)|\,dz \right)\,dr = \omega_n \Waps \left( \tfrac{2^{\alpha - n}}{n - \alpha} f\right)(x)\,.
\end{align*}
The reverse pointwise inequality has been well-studied in the case of $\psi(t)=t^{\frac{1}{p-1}},$ $p>1$, see~\cite[Theorem 6.2]{HaMa72}.
Let us point out that the lower bound giving sharpness in Theorem~\ref{theo:boundedness-psi} is shown for nonnegative and radially decreasing functions precisely as in \cite[Theorem~2.1]{Ci-pot}. 

\subsection{Auxiliary function spaces}

\noindent We say that  a measurable function $f : \Omega\to\rn$ belongs to the Morrey space $L^{q}_{\theta}(\Omega)$ for $q\geq 1$ and $\theta\in[0,n]$, 
if and only if \[\|f\|_{L^{q}_{\theta}(\Omega)}:= \sup\limits_{\substack{
B(x_0,R) \subset\rn\\ x_0\in\Omega}} R^{\frac{\theta-n}{q}} \|f\|_{L^q(B_R\cap \Omega)} <\infty\,.\]
 We say that $f : \Omega \to \rn$ belongs to the Lorentz--Morrey space $\Lambda^{t, q}_{\theta}(\Omega)$ for $q\in[1,\infty]$ and $\theta\in[0,n]$,
if and only if\[\|f\|_{\Lambda^{t,q}_{\theta}(\Omega)}:=\sup\limits_{\substack{
B(x_0,R) \subset\rn\\ x_0\in\Omega}}  R^\frac{\theta-n}{t}\|f\|_{\Lambda^{t, q}(B_R\cap \Omega)}<\infty\,.\]
%
We shall consider data in the Orlicz space $L\log L$, where the modular function $t\mapsto t\log(e+t)$ satisfies $\Delta_2$-condition, but is growing essentially less rapidly than $t^{1+\ve}$ for any $\ve>0$. 

\noindent We define the space $L\log  L(\Omega)$ as a subset of integrable functions $f:\Omega\to\R$ such that 
\[\begin{split}\|f\|_{L \log  L(\Omega)}&=\inf\left\{\lambda>0:\quad\int_\Omega\left|\frac{f}{\lambda}\right|\log\left(e +\left|\frac{f}{\lambda}\right| \right)dx\leq 1\right\} <\infty\,.\end{split}\]

The classical reference for boundedness of maximal function is~\cite{adams-hedberg}, most of the needed estimates can be found in~\cite{min-grad-est}. We use the following estimates on maximal operators. 
 

\begin{lem}\label{lem:maximal-function} 
 Let $B\subset \rn$ be a ball, $n\geq 2$, and  $f:\rn\to\R$  is a locally integrable function supported in $B$. Then
 \begin{itemize}
 \item[i)]
If $1<q<n$ and $s \in(0,\infty]$, then there exists $c = c(n,q,s)$ such that for every it holds that \[\| M_1 f\|_{\Lambda^{\frac{n q}{n-q},s}(B)}\leq\ c \|f \|_{\Lambda^{q,s}(B)}\,.\] 
 \item[ii)]
There exists $c = c(n)$ such that for every it holds that
 \[\| M_1 f\|_{L^\frac{n}{n-1}(B)}\leq c(n)|B|^\frac{1}{n} \| f\|_{L \log L(B)}\,.\]
 \item[iii)] If $1<q<\theta\leq n$,  then there exists $c >0$ such that for every it holds that
\[\| M_1 f\|_{L^{\frac{\theta q}{\theta-q}}_{\theta}(B)}\leq c \|f\|_{L^{q}_{\theta}(B)}\,.\]
 \item[iv)] If $1<q<\theta\leq n$ and $s\in(0,\infty)$,  then there exists $c >0$ such that for every it holds that
\[\|M_1 f\|_{\Lambda^{\frac{\theta q}{\theta-q},\frac{\theta s}{\theta-q}}(B)}\leq c \|f\|_{\Lambda^{q, s}_{\theta}(B)}\,.\]
\end{itemize} 
\end{lem}

\subsection{Technical lemmas}
In this section we assume that $G$ and $g$ are as in Section~\ref{sec:appl}. 

\begin{lem}\label{lem:psi-rearr}
Suppose that $\psi : \rp \to \rp$ is a non-decreasing and left-continuous function. Then for any function $f$ satisfying~\eqref{skoncz-wart} it holds that
\begin{equation*}
    \left( \psi(|f|)\right)^{*}(t) = \psi(f^{*}(t))\,.
\end{equation*}
\end{lem}
\begin{proof}
We observe that the function $\psi \circ f^*$ is non-increasing and right-continuous. Therefore, it suffices to prove that $\psi \circ f^*$ is equimeasurable with $\psi \circ |f|$. Let us take arbitrary $t > 0$. If $t < \psi(0)$ or $t \geq \lim_{s \to \infty} \psi(s)$, then trivially $|\{x : \psi(|f(x)|) > t\}| = |\{x : \psi(f^*(x)) > t\}|$. In the opposite case, we take the largest number $\tau \leq t$ such that there exist $s$ satisfying $\psi(s) = \tau$. Let us assume that $s$ is the largest number such that $\psi(s) = \tau$. Then we conclude by noting that $    |\{x : \psi(|f(x)|) > t\}| 
    = |\{x : |f(x)| > s\}|
    = |\{x : f^*(x) > s\}| = |\{x : \psi(f^*(x)) > \tau\}| = |\{x : \psi(f^*(x)) > t\}|\,.$
\end{proof}

\begin{lem}\label{lem:gpow}
For any $\beta > s_G - 1$ and any constant $C \geq 0$ there exist $c_{G, 1} = c_{G, 1}(\beta, i_G) > 0$ and $c_{G, 2} = c_{G, 2}(\beta, s_G) > 0$ such that
\begin{equation}\label{eq:gpow}
c_{G, 1} tg^{-1}(Ct^{-\beta}) \leq \int_{t}^{\infty} g^{-1}(Cs^{-\beta})\,ds \leq c_{G, 2} tg^{-1}(Ct^{-\beta})\,.  \end{equation}
\end{lem}

\begin{proof}
We use the fact that function $t \mapsto \frac{g^{-1}(Ct)}{t^{\frac{1}{i_G-1}}}$ is comparable to non-increasing and function $t \mapsto \frac{g^{-1}(Ct)}{t^{\frac{1}{s_G-1}}}$ is comparable to non-decreasing to estimate
\[ \int_{t}^{\infty} g^{-1}(Cs^{-\beta})\,ds = \int_{t}^{\infty} s^{-\beta \frac{1}{i_G-1}}\frac{g^{-1}(Cs^{-\beta})}{s^{-\beta \frac{1}{i_G-1}}}\,ds \geq \frac{g^{-1}(Ct^{-\beta})}{t^{-\beta \frac{1}{i_G-1}}} \int_{t}^{\infty} s^{-\beta \frac{1}{i_G-1}}\,ds = \frac{i_G - 1}{\beta - i_G + 1} tg^{-1}(Ct^{-\beta}) \]
and analogously
\[ \int_{t}^{\infty} g^{-1}(Cs^{-\beta})\,ds = \int_{t}^{\infty} s^{-\beta \frac{1}{s_G-1}}\frac{g^{-1}(Cs^{-\beta})}{s^{-\beta \frac{1}{s_G-1}}} \leq \frac{g^{-1}(Ct^{-\beta})}{t^{-\beta \frac{1}{s_G-1}}}\,ds \int_{t}^{\infty} s^{-\beta \frac{1}{s_G-1}}\,ds = \frac{s_G - 1}{\beta - s_G + 1} tg^{-1}(Ct^{-\beta})\,. \]
These two inequalities give us~\eqref{eq:gpow} with $c_{G, 1} = \tfrac{i_G - 1}{\beta - i_G + 1}$ and $c_{G, 2} = \tfrac{s_G - 1}{\beta - s_G + 1}$.
\end{proof}

\begin{lem}\label{lem:jensen}
Suppose that $h:\rp\to\rp$ is a nonnegative function such that the function $r \mapsto r^{\gamma}h(r)$ is convex for some $\gamma \geq 0$. Then for any non-increasing, nonnegative function $\phi$ the following inequality holds true
\begin{equation}\label{eq:jensen} h\left( \frac{1}{t}\int_{0}^{t}\phi(s)ds \right) \leq t^{\gamma - 1}\int_{0}^{t}s^{-\gamma}h(\phi(s))ds\,.  \end{equation}
\end{lem}

\begin{proof}
By Jensen's inequality for function $r \mapsto r^{\gamma}h(r)$ we have
\[ \left(\frac{1}{t}\int_{0}^{t}\phi(s)\,ds \right)^{\gamma}h\left( \frac{1}{t}\int_{0}^{t}\phi(s)\,ds\right) \leq \frac{1}{t}\int_{0}^{t}(\phi(s))^{\gamma}h(\phi(s))\,ds\,, \]
which gives us
\begin{align*}
     h\left( \frac{1}{t}\int_{0}^{t}\phi(s)\,ds\right) 
     &\leq t^{\gamma - 1}\int_{0}^{t} \left( \frac{\phi(s)}{\frac{1}{s}\int_{0}^{s}\phi(r)\,dr} \right)^{\gamma}s^{-\gamma}h(\phi(s))\,ds \leq t^{\gamma - 1}\int_{0}^{t} s^{-\gamma}h(\phi(s))\,ds\,. 
\end{align*}

\end{proof}
\begin{lem}\label{lem:gconv} The function $t \mapsto t^{1 - \frac{1}{s_G - 1}}g^{-1}(t)$ is comparable to a convex function.
\end{lem}
\begin{proof}
Since the function $t \mapsto h_{s_G}(t):=t^{-\frac{1}{s_G-1}}g^{-1}(t)$ is comparable to a non-decreasing function, the function $t\mapsto\int_0^t h_{s_G}(s)\,ds$
%
%
is comparable to a convex function. Observe that
$\int_{0}^{t}h_{s_G}(s)\,ds \leq t h_{s_G}(t)\,.$ 
Since $t \mapsto h_{i_G}(t):= t^{-\frac{1}{i_G-1}}g^{-1}(t)$ is comparable to a non-increasing function, we have
\[ \int_{0}^{t} h_{s_G}(s)\,ds = \int_{0}^{t} s^{\frac{1}{i_G - 1} - \frac{1}{s_G - 1}} s^{-\frac{1}{i_G-1}}g^{-1}(s)\,ds \geq ct h_{s_G}(t)\,,\] 
where in the last inequality we used the fact that $\int_{0}^{t} s^{\frac{1}{i_G-1} - \frac{1}{s_G-1}}\,ds < \infty$, which holds true by~\eqref{iG-sG}. 
Therefore which means that $t \mapsto t h_{s_G}(t)$ is comparable to a convex function (which is comparable to $t\mapsto\int_0^t h_{s_G}(s)\,ds$).
\end{proof}
\subsection{Inverses}

\begin{lem}\label{lem:zyginv}
If $1 < p < \infty$, $\alpha \in \R$, then there exist $s_0 > 1$ such that for $s \geq s_0$ we have that function $G_s(t)=t^p(\log(s+t))^\alpha,$ $1<p<\infty,\ \alpha\in\R$ satisfies $G_s \in \Delta_2 \cap \nabla_2$ and
\begin{equation*} g_s(t) \approx t^{p-1}\log^{\alpha}{(s + t)} \qquad \text{and}\qquad g_s^{-1}(t)\approx t^{\frac{1}{p - 1}}(\log{(s + t)})^{-\frac{\alpha}{p - 1}}\,.\end{equation*}
Moreover, for all $s_1, s_2 \geq s_0$ we have $G_{s_1} \approx G_{s_2}$, $g_{s_1} \approx g_{s_2}$. As $s \to \infty$, indices $s_{G_s} \geq p$ and $i_{G_{s}} \leq p$ can be arbitrarily close to $p$.
\end{lem}

\begin{proof} 
Let us notice that
\begin{equation*}
    g_{s}(t) = t^{p - 1}\log^{\alpha}(s + t)\left(\frac{\alpha t}{(s + t)\log{(s+t)}} + p\right)\,,
\end{equation*}
which means that
\begin{equation*}
    \frac{tg_{s}(t)}{G_s(t)} = p + \frac{\alpha t}{(s+t)\log(s+t)}\,.
\end{equation*}
Since $\left| \frac{\alpha t}{(s+t)\log(s+t)} \right| \leq \frac{|\alpha|}{\log(s+t)} \xrightarrow[]{s \to \infty} 0$, for arbitrary $\epsilon > 0$ and for sufficiently large $s$, we have $\frac{tg_{s}(t)}{\overline{G}_s(t)} \in (p - \epsilon, p + \epsilon)$. In particular, for sufficiently large $s$ we have $G_s \in \Delta_2 \cap \nabla_2$. Observe that
\begin{equation}
    \frac{g_s\left( t^{\frac{1}{p - 1}}(\log{(s + t)})^{-\frac{\alpha}{p - 1}}\right)}{t} \approx \left( \frac{\log \left(s + t^{\frac{1}{p-1}} \left( \log(s + t) \right)^{-\frac{\alpha}{p-1}} \right)}{\log(s + t)} \right)^{\alpha} =: l_s(t)\,.
\end{equation}
Since $l_s > 0$, $l_s(0) = 1$ and $l_s$ has finite, positive limit at infinity, we have that $l_s \approx 1$ and by doubling properties of $g_s$ we have $    g_s^{-1}(t) \approx t^{\frac{1}{p-1}}\left(\log (s + t) \right)^{-\frac{\alpha}{p-1}}$.
\end{proof}

\begin{lem}\label{lem:zyginv2}
If $1 < p < \infty$, $\alpha \in \R$, then there exist $s_0 > e$ such that for $s \geq s_0$ we have that function $G_s(t)=t^p(\log \log(s+t))^\alpha,$ $1<p<\infty,\ \alpha\in\R$ satisfies $G_s \in \Delta_2 \cap \nabla_2$ and
\begin{equation*} g_s(t) \approx t^{p-1}(\log \log{(s + t)})^{\alpha} \qquad \text{and}\qquad g_s^{-1}(t)\approx t^{\frac{1}{p - 1}}(\log \log{(s + t)})^{-\frac{\alpha}{p - 1}}\,.\end{equation*}
Moreover, for all $s_1, s_2 \geq s_0$ we have $G_{s_1} \approx G_{s_2}$, $g_{s_1} \approx g_{s_2}$. As $s \to \infty$, indices $s_{G_s} \geq p$ and $i_{G_{s}} \leq p$ can be arbitrarily close to $p$.
\end{lem}

\begin{proof} 
Let us notice that
\begin{equation*}
    g_{s}(t) = t^{p - 1}(\log \log(s + t))^{\alpha}\left(\frac{\alpha t}{(s + t)\log{(s+t)}\log \log (s+t)} + p\right)\,,
\end{equation*}
which means that
\begin{equation*}
    \frac{tg_{s}(t)}{G_s(t)} = p + \frac{\alpha t}{(s+t)\log(s+t)\log \log (s+t)}\,.
\end{equation*}
Since $\left| \frac{\alpha t}{(s+t)\log(s+t)\log \log (s+t)} \right| \leq \frac{|\alpha|}{\log(s+t)\log \log (s+t)} \xrightarrow[]{s \to \infty} 0$, for arbitrary $\epsilon > 0$ and for sufficiently large $s$, we have $\frac{tg_{s}(t)}{\overline{G}_s(t)} \in (p - \epsilon, p + \epsilon)$. In particular, for sufficiently large $s$ we have $G_s \in \Delta_2 \cap \nabla_2$. Observe that
\begin{equation}
    \frac{g_s\left( t^{\frac{1}{p - 1}}(\log \log{(s + t)})^{-\frac{\alpha}{p - 1}}\right)}{t} \approx \left( \frac{\log \log \left(s + t^{\frac{1}{p-1}} \left( \log \log (s + t) \right)^{-\frac{\alpha}{p-1}} \right)}{\log \log(s + t)} \right)^{\alpha} =: l_s(t)\,.
\end{equation}
Since $l_s > 0$, $l_s(0) = 1$ and $l_s$ has finite, positive limit at infinity, we have that $l_s \approx 1$ and by doubling properties of $g_s$ we have $ g_s^{-1}(t) \approx t^{\frac{1}{p-1}}\left(\log \log (s + t) \right)^{-\frac{\alpha}{p-1}}$.
\end{proof}

\bibliographystyle{plain}
\bibliography{bib}

\end{document}